\documentclass[12pt]{amsart}
\usepackage{amssymb}

\newlength{\horzstretch}
\newlength{\vertstretch}
\theoremstyle{plain}
\newtheorem{lemma}{Lemma}[section]

\newtheorem{theorem}[lemma]{Theorem}
\newtheorem{corollary}[lemma]{Corollary}

\theoremstyle{definition}
\newtheorem{definition}[lemma]{Definition}

\newtheorem{remark}[lemma]{Remark}

\newcommand{\ggcaffil}[1]{\dedicatory{\textup{\larger{#1}}}}

\newcommand{\ggcqedsymbol}{$\square$}
\newcommand{\ggcqed}{\hbox{}\nobreak\hbox{\quad\ggcqedsymbol}}
\newcommand{\ggcendpf}{\ggcqed}
\newcommand{\ggcnopf}{\ggcqed}

\newcommand{\ggcenddef}{\ggcqed}

\newcommand{\ggcendrem}{\ggcenddef}

\newcommand{\setmin}{\setminus}
\newcommand{\lcover}{\mathbin{<\!\!\cdot}}
\newcommand{\cC}{\ensuremath{\mathcal{C}}}

\newcommand{\ggcoverline}[1]{\mskip2mu\overline{\mskip-2mu#1}}
\newcommand{\oG}{\ggcoverline{G}}

\begin{document}
\title[Polyunsaturated Posets and Graphs]
  {Polyunsaturated Posets and Graphs\\
  and the Greene-Kleitman Theorem}
\author{Glenn G. Chappell}
\ggcaffil{Department of Mathematics, Southeast Missouri
  State University}
\address{Department of Mathematics\\
  Southeast Missouri State University\\
  Cape\break Girardeau, MO 63701\\
  USA}
\email{gchappell@semovm.semo.edu}
\thanks{This work first appeared in the author's
  Ph.D. thesis~\cite{ChaG96},
  supervised by Douglas B. West at the University of Illinois;
  research supported by U.S. Department of Education
  grant DE-P200A40311.}
\date{July 31, 1998}
\begin{abstract}
  A partition of a finite poset into chains places a natural upper bound
  on the size of a union of $k$ antichains.
  A chain partition is
  \emph{$k$-saturated} if this bound is achieved.
  Greene and Kleitman~\cite{GrKl76} proved that, for each $k$,
  every finite poset has a simultaneously
  $k$- and $k+1$-saturated chain partition.
  West~\cite{WesD86} showed that the Greene-Kleitman Theorem is
  best-possible in a strong sense by exhibiting, for each $c\ge4$,
  a poset with longest chain of cardinality $c$ and no
  $k$- and $l$-saturated chain partition for any distinct, nonconsecutive
  $k,l<c$.
  We call such posets \emph{polyunsaturated}.
  We give necessary and sufficient conditions for the existence of
  polyunsaturated posets with prescribed height, width, and cardinality.
  We prove these results in the more general context of graphs satisfying
  an analogue of the Greene-Kleitman Theorem.
  Lastly, we discuss analogous results for antichain partitions.
\end{abstract}
\maketitle

\section{Introduction} \label{S:intro}

Let $P$ be a finite partially ordered set.
A theorem of Dilworth~\cite{DilR50} states that
the cardinality of the largest antichain in $P$ is equal to the
least number of chains into which $P$ can be partitioned.

Greene and Kleitman~\cite{GrKl76} generalized this to unions
of antichains.
A \emph{$k$-family} in $P$
is a union of $k$ antichains;
we denote the size of the largest $k$-family by $d_k(P)$.
Given a partition $\cC$ of a set $S$,
we define the \emph{$k$-norm} of $\cC$, denoted
$m_k(\cC)$, as follows:
\[
m_k(\cC):=\sum_{C\in{\cC}}\min\left\{k,|C|\right\}.
\]
For each chain partition $\cC$ of $P$,
$m_k(\cC)$ is an upper bound for $d_k(P)$;
a chain partition $\cC$ is \emph{$k$-saturated} if
this bound is achieved.
Dilworth's Theorem says that
every finite poset $P$ has a $1$-saturated chain partition.
Greene and Kleitman~\cite[Thm.~3.11]{GrKl76}
proved that for each $k$,
$P$ has a simultaneously $k$- and
$k+1$-saturated chain partition.
Other proofs of this result can be found
in~\cite{CaEd92,FomS78,FraA80,HoSc77,PerH84,SakM79}.

\begin{figure}[htbp]
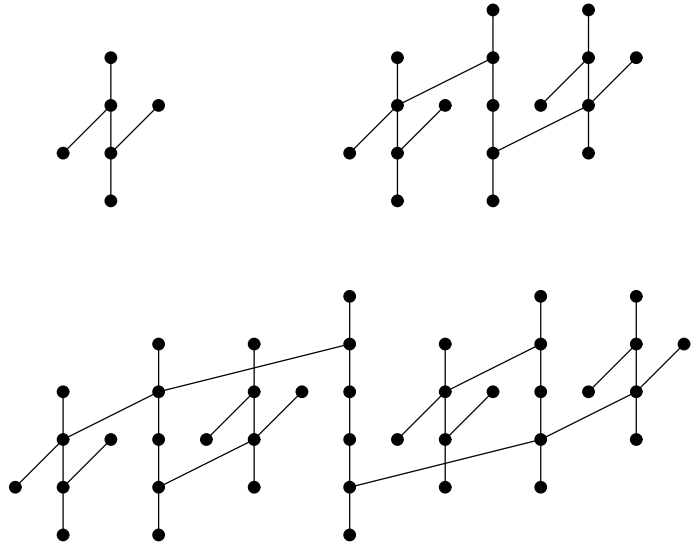

\expandafter\ifx\csname graph\endcsname\relax \csname newbox\endcsname\graph\fi
\expandafter\ifx\csname graphtemp\endcsname\relax \csname newdimen\endcsname\graphtemp\fi
\setbox\graph=\vtop{\vskip 0pt\hbox{%
\fontsize{10}{12}\selectfont 
    \special{pn 8}%
    \special{sh 1.000}%
    \special{ar 281 781 31 31 0 6.28319}%
    \special{sh 1.000}%
    \special{ar 531 1031 31 31 0 6.28319}%
    \special{sh 1.000}%
    \special{ar 531 781 31 31 0 6.28319}%
    \special{sh 1.000}%
    \special{ar 531 531 31 31 0 6.28319}%
    \special{sh 1.000}%
    \special{ar 531 281 31 31 0 6.28319}%
    \special{sh 1.000}%
    \special{ar 781 531 31 31 0 6.28319}%
    \special{pa 281 781}%
    \special{pa 531 531}%
    \special{fp}%
    \special{pa 531 1031}%
    \special{pa 531 281}%
    \special{fp}%
    \special{pa 531 781}%
    \special{pa 781 531}%
    \special{fp}%
    \special{sh 1.000}%
    \special{ar 1781 781 31 31 0 6.28319}%
    \special{sh 1.000}%
    \special{ar 2031 1031 31 31 0 6.28319}%
    \special{sh 1.000}%
    \special{ar 2031 781 31 31 0 6.28319}%
    \special{sh 1.000}%
    \special{ar 2031 531 31 31 0 6.28319}%
    \special{sh 1.000}%
    \special{ar 2031 281 31 31 0 6.28319}%
    \special{sh 1.000}%
    \special{ar 2281 531 31 31 0 6.28319}%
    \special{sh 1.000}%
    \special{ar 2531 1031 31 31 0 6.28319}%
    \special{sh 1.000}%
    \special{ar 2531 781 31 31 0 6.28319}%
    \special{sh 1.000}%
    \special{ar 2531 531 31 31 0 6.28319}%
    \special{sh 1.000}%
    \special{ar 2531 281 31 31 0 6.28319}%
    \special{sh 1.000}%
    \special{ar 2531 31 31 31 0 6.28319}%
    \special{sh 1.000}%
    \special{ar 2781 531 31 31 0 6.28319}%
    \special{sh 1.000}%
    \special{ar 3031 781 31 31 0 6.28319}%
    \special{sh 1.000}%
    \special{ar 3031 531 31 31 0 6.28319}%
    \special{sh 1.000}%
    \special{ar 3031 281 31 31 0 6.28319}%
    \special{sh 1.000}%
    \special{ar 3031 31 31 31 0 6.28319}%
    \special{sh 1.000}%
    \special{ar 3281 281 31 31 0 6.28319}%
    \special{pa 1781 781}%
    \special{pa 2031 531}%
    \special{fp}%
    \special{pa 2031 1031}%
    \special{pa 2031 281}%
    \special{fp}%
    \special{pa 2031 781}%
    \special{pa 2281 531}%
    \special{fp}%
    \special{pa 2031 531}%
    \special{pa 2531 281}%
    \special{fp}%
    \special{pa 2531 1031}%
    \special{pa 2531 31}%
    \special{fp}%
    \special{pa 2531 781}%
    \special{pa 3031 531}%
    \special{fp}%
    \special{pa 2781 531}%
    \special{pa 3031 281}%
    \special{fp}%
    \special{pa 3031 781}%
    \special{pa 3031 31}%
    \special{fp}%
    \special{pa 3031 531}%
    \special{pa 3281 281}%
    \special{fp}%
    \special{sh 1.000}%
    \special{ar 31 2531 31 31 0 6.28319}%
    \special{sh 1.000}%
    \special{ar 281 2781 31 31 0 6.28319}%
    \special{sh 1.000}%
    \special{ar 281 2531 31 31 0 6.28319}%
    \special{sh 1.000}%
    \special{ar 281 2281 31 31 0 6.28319}%
    \special{sh 1.000}%
    \special{ar 281 2031 31 31 0 6.28319}%
    \special{sh 1.000}%
    \special{ar 531 2281 31 31 0 6.28319}%
    \special{sh 1.000}%
    \special{ar 781 2781 31 31 0 6.28319}%
    \special{sh 1.000}%
    \special{ar 781 2531 31 31 0 6.28319}%
    \special{sh 1.000}%
    \special{ar 781 2281 31 31 0 6.28319}%
    \special{sh 1.000}%
    \special{ar 781 2031 31 31 0 6.28319}%
    \special{sh 1.000}%
    \special{ar 781 1781 31 31 0 6.28319}%
    \special{sh 1.000}%
    \special{ar 1031 2281 31 31 0 6.28319}%
    \special{sh 1.000}%
    \special{ar 1281 2531 31 31 0 6.28319}%
    \special{sh 1.000}%
    \special{ar 1281 2281 31 31 0 6.28319}%
    \special{sh 1.000}%
    \special{ar 1281 2031 31 31 0 6.28319}%
    \special{sh 1.000}%
    \special{ar 1281 1781 31 31 0 6.28319}%
    \special{sh 1.000}%
    \special{ar 1531 2031 31 31 0 6.28319}%
    \special{sh 1.000}%
    \special{ar 1781 2781 31 31 0 6.28319}%
    \special{sh 1.000}%
    \special{ar 1781 2531 31 31 0 6.28319}%
    \special{sh 1.000}%
    \special{ar 1781 2281 31 31 0 6.28319}%
    \special{sh 1.000}%
    \special{ar 1781 2031 31 31 0 6.28319}%
    \special{sh 1.000}%
    \special{ar 1781 1781 31 31 0 6.28319}%
    \special{sh 1.000}%
    \special{ar 1781 1531 31 31 0 6.28319}%
    \special{sh 1.000}%
    \special{ar 2031 2281 31 31 0 6.28319}%
    \special{sh 1.000}%
    \special{ar 2281 2531 31 31 0 6.28319}%
    \special{sh 1.000}%
    \special{ar 2281 2281 31 31 0 6.28319}%
    \special{sh 1.000}%
    \special{ar 2281 2031 31 31 0 6.28319}%
    \special{sh 1.000}%
    \special{ar 2281 1781 31 31 0 6.28319}%
    \special{sh 1.000}%
    \special{ar 2531 2031 31 31 0 6.28319}%
    \special{sh 1.000}%
    \special{ar 2781 2531 31 31 0 6.28319}%
    \special{sh 1.000}%
    \special{ar 2781 2281 31 31 0 6.28319}%
    \special{sh 1.000}%
    \special{ar 2781 2031 31 31 0 6.28319}%
    \special{sh 1.000}%
    \special{ar 2781 1781 31 31 0 6.28319}%
    \special{sh 1.000}%
    \special{ar 2781 1531 31 31 0 6.28319}%
    \special{sh 1.000}%
    \special{ar 3031 2031 31 31 0 6.28319}%
    \special{sh 1.000}%
    \special{ar 3281 2281 31 31 0 6.28319}%
    \special{sh 1.000}%
    \special{ar 3281 2031 31 31 0 6.28319}%
    \special{sh 1.000}%
    \special{ar 3281 1781 31 31 0 6.28319}%
    \special{sh 1.000}%
    \special{ar 3281 1531 31 31 0 6.28319}%
    \special{sh 1.000}%
    \special{ar 3531 1781 31 31 0 6.28319}%
    \special{pa 31 2531}%
    \special{pa 281 2281}%
    \special{fp}%
    \special{pa 281 2781}%
    \special{pa 281 2031}%
    \special{fp}%
    \special{pa 281 2531}%
    \special{pa 531 2281}%
    \special{fp}%
    \special{pa 281 2281}%
    \special{pa 781 2031}%
    \special{fp}%
    \special{pa 781 2781}%
    \special{pa 781 1781}%
    \special{fp}%
    \special{pa 781 2531}%
    \special{pa 1281 2281}%
    \special{fp}%
    \special{pa 1031 2281}%
    \special{pa 1281 2031}%
    \special{fp}%
    \special{pa 1281 2531}%
    \special{pa 1281 1781}%
    \special{fp}%
    \special{pa 1281 2281}%
    \special{pa 1531 2031}%
    \special{fp}%
    \special{pa 781 2031}%
    \special{pa 1781 1781}%
    \special{fp}%
    \special{pa 1781 2781}%
    \special{pa 1781 1531}%
    \special{fp}%
    \special{pa 1781 2531}%
    \special{pa 2781 2281}%
    \special{fp}%
    \special{pa 2031 2281}%
    \special{pa 2281 2031}%
    \special{fp}%
    \special{pa 2281 2531}%
    \special{pa 2281 1781}%
    \special{fp}%
    \special{pa 2281 2281}%
    \special{pa 2531 2031}%
    \special{fp}%
    \special{pa 2281 2031}%
    \special{pa 2781 1781}%
    \special{fp}%
    \special{pa 2781 2531}%
    \special{pa 2781 1531}%
    \special{fp}%
    \special{pa 2781 2281}%
    \special{pa 3281 2031}%
    \special{fp}%
    \special{pa 3031 2031}%
    \special{pa 3281 1781}%
    \special{fp}%
    \special{pa 3281 2281}%
    \special{pa 3281 1531}%
    \special{fp}%
    \special{pa 3281 2031}%
    \special{pa 3531 1781}%
    \special{fp}%
    \hbox{\vrule depth2.813in width0pt height 0pt}%
    \kern 3.563in
  }%
}%
\ht\graph=\dp\graph
\dp\graph=0pt

\[\box\graph\]
\caption{West's posets.}
\label{F:wests}
\end{figure}

It is natural to ask whether a stronger result is possible:
can we always find a $k$-,
$k+1$-, and $k+2$-saturated chain partition?
Greene and Kleitman showed that the answer to this question is
``no'' by exhibiting a poset
(the first poset in Figure~\ref{F:wests})
with no $1$- and $3$-saturated partition.

The \emph{height} $h(P)$ of a poset $P$ is the maximum size of
a chain in $P$.
West~\cite{WesD86} generalized the above example by
exhibiting a poset of each
height $c\ge4$ in which no chain partition is
$k$- and $l$-saturated for any distinct, nonconsecutive $k,l<c$
(Figure~\ref{F:wests}).
We say such a poset is \emph{polyunsaturated}.
Since every chain partition of $P$ is $k$-saturated
for each $k\ge h(P)$, West's examples show that
the Greene-Kleitman Theorem is best-possible in a strong sense.

\begin{figure}[htbp]
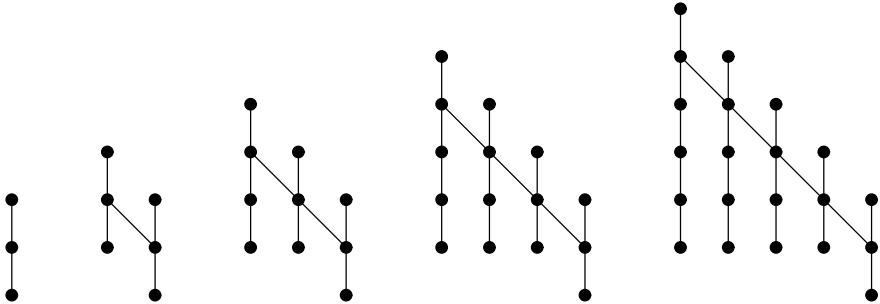

\expandafter\ifx\csname graph\endcsname\relax \csname newbox\endcsname\graph\fi
\expandafter\ifx\csname graphtemp\endcsname\relax \csname newdimen\endcsname\graphtemp\fi
\setbox\graph=\vtop{\vskip 0pt\hbox{%
\fontsize{10}{12}\selectfont 
    \special{pn 8}%
    \special{sh 1.000}%
    \special{ar 31 1531 31 31 0 6.28319}%
    \special{sh 1.000}%
    \special{ar 31 1281 31 31 0 6.28319}%
    \special{sh 1.000}%
    \special{ar 31 1031 31 31 0 6.28319}%
    \special{pa 31 1531}%
    \special{pa 31 1031}%
    \special{fp}%
    \special{sh 1.000}%
    \special{ar 531 1281 31 31 0 6.28319}%
    \special{sh 1.000}%
    \special{ar 531 1031 31 31 0 6.28319}%
    \special{sh 1.000}%
    \special{ar 531 781 31 31 0 6.28319}%
    \special{sh 1.000}%
    \special{ar 781 1531 31 31 0 6.28319}%
    \special{sh 1.000}%
    \special{ar 781 1281 31 31 0 6.28319}%
    \special{sh 1.000}%
    \special{ar 781 1031 31 31 0 6.28319}%
    \special{pa 531 1281}%
    \special{pa 531 781}%
    \special{fp}%
    \special{pa 781 1531}%
    \special{pa 781 1031}%
    \special{fp}%
    \special{pa 531 1031}%
    \special{pa 781 1281}%
    \special{fp}%
    \special{sh 1.000}%
    \special{ar 1281 1281 31 31 0 6.28319}%
    \special{sh 1.000}%
    \special{ar 1281 1031 31 31 0 6.28319}%
    \special{sh 1.000}%
    \special{ar 1281 781 31 31 0 6.28319}%
    \special{sh 1.000}%
    \special{ar 1281 531 31 31 0 6.28319}%
    \special{sh 1.000}%
    \special{ar 1531 1281 31 31 0 6.28319}%
    \special{sh 1.000}%
    \special{ar 1531 1031 31 31 0 6.28319}%
    \special{sh 1.000}%
    \special{ar 1531 781 31 31 0 6.28319}%
    \special{sh 1.000}%
    \special{ar 1781 1531 31 31 0 6.28319}%
    \special{sh 1.000}%
    \special{ar 1781 1281 31 31 0 6.28319}%
    \special{sh 1.000}%
    \special{ar 1781 1031 31 31 0 6.28319}%
    \special{pa 1281 1281}%
    \special{pa 1281 531}%
    \special{fp}%
    \special{pa 1531 1281}%
    \special{pa 1531 781}%
    \special{fp}%
    \special{pa 1781 1531}%
    \special{pa 1781 1031}%
    \special{fp}%
    \special{pa 1281 781}%
    \special{pa 1781 1281}%
    \special{fp}%
    \special{sh 1.000}%
    \special{ar 2281 1281 31 31 0 6.28319}%
    \special{sh 1.000}%
    \special{ar 2281 1031 31 31 0 6.28319}%
    \special{sh 1.000}%
    \special{ar 2281 781 31 31 0 6.28319}%
    \special{sh 1.000}%
    \special{ar 2281 531 31 31 0 6.28319}%
    \special{sh 1.000}%
    \special{ar 2281 281 31 31 0 6.28319}%
    \special{sh 1.000}%
    \special{ar 2531 1281 31 31 0 6.28319}%
    \special{sh 1.000}%
    \special{ar 2531 1031 31 31 0 6.28319}%
    \special{sh 1.000}%
    \special{ar 2531 781 31 31 0 6.28319}%
    \special{sh 1.000}%
    \special{ar 2531 531 31 31 0 6.28319}%
    \special{sh 1.000}%
    \special{ar 2781 1281 31 31 0 6.28319}%
    \special{sh 1.000}%
    \special{ar 2781 1031 31 31 0 6.28319}%
    \special{sh 1.000}%
    \special{ar 2781 781 31 31 0 6.28319}%
    \special{sh 1.000}%
    \special{ar 3031 1531 31 31 0 6.28319}%
    \special{sh 1.000}%
    \special{ar 3031 1281 31 31 0 6.28319}%
    \special{sh 1.000}%
    \special{ar 3031 1031 31 31 0 6.28319}%
    \special{pa 2281 1281}%
    \special{pa 2281 281}%
    \special{fp}%
    \special{pa 2531 1281}%
    \special{pa 2531 531}%
    \special{fp}%
    \special{pa 2781 1281}%
    \special{pa 2781 781}%
    \special{fp}%
    \special{pa 3031 1531}%
    \special{pa 3031 1031}%
    \special{fp}%
    \special{pa 2281 531}%
    \special{pa 3031 1281}%
    \special{fp}%
    \special{sh 1.000}%
    \special{ar 3531 1281 31 31 0 6.28319}%
    \special{sh 1.000}%
    \special{ar 3531 1031 31 31 0 6.28319}%
    \special{sh 1.000}%
    \special{ar 3531 781 31 31 0 6.28319}%
    \special{sh 1.000}%
    \special{ar 3531 531 31 31 0 6.28319}%
    \special{sh 1.000}%
    \special{ar 3531 281 31 31 0 6.28319}%
    \special{sh 1.000}%
    \special{ar 3531 31 31 31 0 6.28319}%
    \special{sh 1.000}%
    \special{ar 3781 1281 31 31 0 6.28319}%
    \special{sh 1.000}%
    \special{ar 3781 1031 31 31 0 6.28319}%
    \special{sh 1.000}%
    \special{ar 3781 781 31 31 0 6.28319}%
    \special{sh 1.000}%
    \special{ar 3781 531 31 31 0 6.28319}%
    \special{sh 1.000}%
    \special{ar 3781 281 31 31 0 6.28319}%
    \special{sh 1.000}%
    \special{ar 4031 1281 31 31 0 6.28319}%
    \special{sh 1.000}%
    \special{ar 4031 1031 31 31 0 6.28319}%
    \special{sh 1.000}%
    \special{ar 4031 781 31 31 0 6.28319}%
    \special{sh 1.000}%
    \special{ar 4031 531 31 31 0 6.28319}%
    \special{sh 1.000}%
    \special{ar 4281 1281 31 31 0 6.28319}%
    \special{sh 1.000}%
    \special{ar 4281 1031 31 31 0 6.28319}%
    \special{sh 1.000}%
    \special{ar 4281 781 31 31 0 6.28319}%
    \special{sh 1.000}%
    \special{ar 4531 1531 31 31 0 6.28319}%
    \special{sh 1.000}%
    \special{ar 4531 1281 31 31 0 6.28319}%
    \special{sh 1.000}%
    \special{ar 4531 1031 31 31 0 6.28319}%
    \special{pa 3531 1281}%
    \special{pa 3531 31}%
    \special{fp}%
    \special{pa 3781 1281}%
    \special{pa 3781 281}%
    \special{fp}%
    \special{pa 4031 1281}%
    \special{pa 4031 531}%
    \special{fp}%
    \special{pa 4281 1281}%
    \special{pa 4281 781}%
    \special{fp}%
    \special{pa 4531 1531}%
    \special{pa 4531 1031}%
    \special{fp}%
    \special{pa 3531 281}%
    \special{pa 4531 1281}%
    \special{fp}%
    \hbox{\vrule depth1.563in width0pt height 0pt}%
    \kern 4.563in
  }%
}%
\ht\graph=\dp\graph
\dp\graph=0pt

\[\box\graph\]
\caption{Our posets.}
\label{F:mine}
\end{figure}

West's examples have exponential width
as a function of their heights.
In Section~\ref{S:construction}
we construct much narrower polyunsaturated posets
(Figure~\ref{F:mine}).
We show later that these posets have the smallest possible width and
cardinality for polyunsaturated posets with the same height.

We denote the comparability graph of a poset $P$ by $G(P)$.
The chains and antichains of $P$ are precisely the cliques
and independent sets of $G(P)$.
In Section~\ref{S:mainres},
we use these ideas to extend the concept of
polyunsaturation to graphs that satisfy an analogue of the
Greene-Kleitman Theorem, even though they may not be
comparability graphs.
Thus, while our results are primarily of interest in
the context of posets, we prove them in
a more general graph-theoretic form.

Our main result
(Theorem~\ref{T:main})
characterizes, in graph-theoretic form,
the possible values of $d_k(P)$ for a polyunsaturated poset $P$.
In Section~\ref{S:icn}
we obtain corollaries of Theorem~\ref{T:main}
that give necessary and sufficient conditions for the existence
of polyunsaturated posets and graphs with certain parameters.
In Section~\ref{S:indsetpart}
we discuss results analogous to those above for partitions of a poset
into antichains.

\section{The Construction} \label{S:construction}

We now construct the narrowest polyunsaturated posets, relative to height.
When $x<y$ is a cover relation we write $x\lcover y$.
We denote a finite sequence of positive integers by an underlined letter,
and we use subscripted letters to name the elements of the sequence.
Given a sequence $\underline b=\left(b_1,b_2,\dotsc,b_t\right)$, the
\emph{difference sequence} of
$\underline b$ is the sequence $\Delta\underline b$, where
$\Delta b_i=b_i-b_{i-1}$, with the convention that $\Delta b_1=b_1$.

\begin{figure}[htbp]
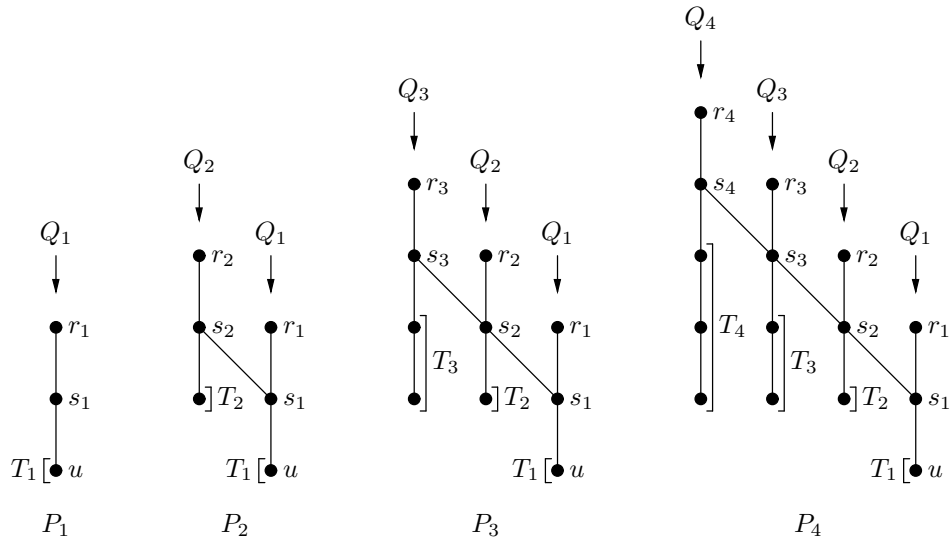

\expandafter\ifx\csname graph\endcsname\relax \csname newbox\endcsname\graph\fi
\expandafter\ifx\csname graphtemp\endcsname\relax \csname newdimen\endcsname\graphtemp\fi
\setbox\graph=\vtop{\vskip 0pt\hbox{%
\fontsize{10}{12}\selectfont 
    \graphtemp=\baselineskip\multiply\graphtemp by -1\divide\graphtemp by 2
    \advance\graphtemp by .5ex\advance\graphtemp by 2.656in
    \rlap{\kern 0.094in\lower\graphtemp\hbox to 0pt{\hss $P_1$\hss}}%
    \special{pn 8}%
    \special{sh 1.000}%
    \special{ar 94 2281 31 31 0 6.28319}%
    \special{sh 1.000}%
    \special{ar 94 1906 31 31 0 6.28319}%
    \special{sh 1.000}%
    \special{ar 94 1531 31 31 0 6.28319}%
    \special{pa 94 2281}%
    \special{pa 94 1531}%
    \special{fp}%
    \graphtemp=\baselineskip\multiply\graphtemp by -1\divide\graphtemp by 2
    \advance\graphtemp by .5ex\advance\graphtemp by 1.125in
    \rlap{\kern 0.094in\lower\graphtemp\hbox to 0pt{\hss $Q_1$\hss}}%
    \special{pa 94 1156}%
    \special{pa 94 1344}%
    \special{fp}%
    \special{sh 1.000}%
    \special{pa 109 1281}%
    \special{pa 94 1344}%
    \special{pa 78 1281}%
    \special{pa 109 1281}%
    \special{fp}%
    \graphtemp=.5ex\advance\graphtemp by 1.531in
    \rlap{\kern 0.156in\lower\graphtemp\hbox to 0pt{$r_1$\hss}}%
    \graphtemp=.5ex\advance\graphtemp by 1.906in
    \rlap{\kern 0.156in\lower\graphtemp\hbox to 0pt{$s_1$\hss}}%
    \graphtemp=.5ex\advance\graphtemp by 2.281in
    \rlap{\kern 0.156in\lower\graphtemp\hbox to 0pt{$u$\hss}}%
    \graphtemp=.5ex\advance\graphtemp by 2.281in
    \rlap{\kern 0.000in\lower\graphtemp\hbox to 0pt{\hss $T_1$}}%
    \special{pa 63 2219}%
    \special{pa 31 2219}%
    \special{fp}%
    \special{pa 31 2219}%
    \special{pa 31 2344}%
    \special{fp}%
    \special{pa 31 2344}%
    \special{pa 63 2344}%
    \special{fp}%
    \graphtemp=\baselineskip\multiply\graphtemp by -1\divide\graphtemp by 2
    \advance\graphtemp by .5ex\advance\graphtemp by 2.656in
    \rlap{\kern 1.031in\lower\graphtemp\hbox to 0pt{\hss $P_2$\hss}}%
    \special{sh 1.000}%
    \special{ar 844 1906 31 31 0 6.28319}%
    \special{sh 1.000}%
    \special{ar 844 1531 31 31 0 6.28319}%
    \special{sh 1.000}%
    \special{ar 844 1156 31 31 0 6.28319}%
    \special{sh 1.000}%
    \special{ar 1219 2281 31 31 0 6.28319}%
    \special{sh 1.000}%
    \special{ar 1219 1906 31 31 0 6.28319}%
    \special{sh 1.000}%
    \special{ar 1219 1531 31 31 0 6.28319}%
    \special{pa 844 1906}%
    \special{pa 844 1156}%
    \special{fp}%
    \special{pa 1219 2281}%
    \special{pa 1219 1531}%
    \special{fp}%
    \special{pa 844 1531}%
    \special{pa 1219 1906}%
    \special{fp}%
    \graphtemp=\baselineskip\multiply\graphtemp by -1\divide\graphtemp by 2
    \advance\graphtemp by .5ex\advance\graphtemp by 1.125in
    \rlap{\kern 1.219in\lower\graphtemp\hbox to 0pt{\hss $Q_1$\hss}}%
    \special{pa 1219 1156}%
    \special{pa 1219 1344}%
    \special{fp}%
    \special{sh 1.000}%
    \special{pa 1234 1281}%
    \special{pa 1219 1344}%
    \special{pa 1203 1281}%
    \special{pa 1234 1281}%
    \special{fp}%
    \graphtemp=.5ex\advance\graphtemp by 1.531in
    \rlap{\kern 1.281in\lower\graphtemp\hbox to 0pt{$r_1$\hss}}%
    \graphtemp=.5ex\advance\graphtemp by 1.906in
    \rlap{\kern 1.281in\lower\graphtemp\hbox to 0pt{$s_1$\hss}}%
    \graphtemp=.5ex\advance\graphtemp by 2.281in
    \rlap{\kern 1.281in\lower\graphtemp\hbox to 0pt{$u$\hss}}%
    \graphtemp=.5ex\advance\graphtemp by 2.281in
    \rlap{\kern 1.125in\lower\graphtemp\hbox to 0pt{\hss $T_1$}}%
    \special{pa 1188 2219}%
    \special{pa 1156 2219}%
    \special{fp}%
    \special{pa 1156 2219}%
    \special{pa 1156 2344}%
    \special{fp}%
    \special{pa 1156 2344}%
    \special{pa 1188 2344}%
    \special{fp}%
    \graphtemp=\baselineskip\multiply\graphtemp by -1\divide\graphtemp by 2
    \advance\graphtemp by .5ex\advance\graphtemp by 0.750in
    \rlap{\kern 0.844in\lower\graphtemp\hbox to 0pt{\hss $Q_2$\hss}}%
    \special{pa 844 781}%
    \special{pa 844 969}%
    \special{fp}%
    \special{sh 1.000}%
    \special{pa 859 906}%
    \special{pa 844 969}%
    \special{pa 828 906}%
    \special{pa 859 906}%
    \special{fp}%
    \graphtemp=.5ex\advance\graphtemp by 1.156in
    \rlap{\kern 0.906in\lower\graphtemp\hbox to 0pt{$r_2$\hss}}%
    \graphtemp=.5ex\advance\graphtemp by 1.531in
    \rlap{\kern 0.906in\lower\graphtemp\hbox to 0pt{$s_2$\hss}}%
    \graphtemp=.5ex\advance\graphtemp by 1.906in
    \rlap{\kern 0.938in\lower\graphtemp\hbox to 0pt{$T_2$\hss}}%
    \special{pa 875 1844}%
    \special{pa 906 1844}%
    \special{fp}%
    \special{pa 906 1844}%
    \special{pa 906 1969}%
    \special{fp}%
    \special{pa 906 1969}%
    \special{pa 875 1969}%
    \special{fp}%
    \graphtemp=\baselineskip\multiply\graphtemp by -1\divide\graphtemp by 2
    \advance\graphtemp by .5ex\advance\graphtemp by 2.656in
    \rlap{\kern 2.344in\lower\graphtemp\hbox to 0pt{\hss $P_3$\hss}}%
    \special{sh 1.000}%
    \special{ar 1969 1906 31 31 0 6.28319}%
    \special{sh 1.000}%
    \special{ar 1969 1531 31 31 0 6.28319}%
    \special{sh 1.000}%
    \special{ar 1969 1156 31 31 0 6.28319}%
    \special{sh 1.000}%
    \special{ar 1969 781 31 31 0 6.28319}%
    \special{sh 1.000}%
    \special{ar 2344 1906 31 31 0 6.28319}%
    \special{sh 1.000}%
    \special{ar 2344 1531 31 31 0 6.28319}%
    \special{sh 1.000}%
    \special{ar 2344 1156 31 31 0 6.28319}%
    \special{sh 1.000}%
    \special{ar 2719 2281 31 31 0 6.28319}%
    \special{sh 1.000}%
    \special{ar 2719 1906 31 31 0 6.28319}%
    \special{sh 1.000}%
    \special{ar 2719 1531 31 31 0 6.28319}%
    \special{pa 1969 1906}%
    \special{pa 1969 781}%
    \special{fp}%
    \special{pa 2344 1906}%
    \special{pa 2344 1156}%
    \special{fp}%
    \special{pa 2719 2281}%
    \special{pa 2719 1531}%
    \special{fp}%
    \special{pa 1969 1156}%
    \special{pa 2719 1906}%
    \special{fp}%
    \graphtemp=\baselineskip\multiply\graphtemp by -1\divide\graphtemp by 2
    \advance\graphtemp by .5ex\advance\graphtemp by 1.125in
    \rlap{\kern 2.719in\lower\graphtemp\hbox to 0pt{\hss $Q_1$\hss}}%
    \special{pa 2719 1156}%
    \special{pa 2719 1344}%
    \special{fp}%
    \special{sh 1.000}%
    \special{pa 2734 1281}%
    \special{pa 2719 1344}%
    \special{pa 2703 1281}%
    \special{pa 2734 1281}%
    \special{fp}%
    \graphtemp=.5ex\advance\graphtemp by 1.531in
    \rlap{\kern 2.781in\lower\graphtemp\hbox to 0pt{$r_1$\hss}}%
    \graphtemp=.5ex\advance\graphtemp by 1.906in
    \rlap{\kern 2.781in\lower\graphtemp\hbox to 0pt{$s_1$\hss}}%
    \graphtemp=.5ex\advance\graphtemp by 2.281in
    \rlap{\kern 2.781in\lower\graphtemp\hbox to 0pt{$u$\hss}}%
    \graphtemp=.5ex\advance\graphtemp by 2.281in
    \rlap{\kern 2.625in\lower\graphtemp\hbox to 0pt{\hss $T_1$}}%
    \special{pa 2688 2219}%
    \special{pa 2656 2219}%
    \special{fp}%
    \special{pa 2656 2219}%
    \special{pa 2656 2344}%
    \special{fp}%
    \special{pa 2656 2344}%
    \special{pa 2688 2344}%
    \special{fp}%
    \graphtemp=\baselineskip\multiply\graphtemp by -1\divide\graphtemp by 2
    \advance\graphtemp by .5ex\advance\graphtemp by 0.750in
    \rlap{\kern 2.344in\lower\graphtemp\hbox to 0pt{\hss $Q_2$\hss}}%
    \special{pa 2344 781}%
    \special{pa 2344 969}%
    \special{fp}%
    \special{sh 1.000}%
    \special{pa 2359 906}%
    \special{pa 2344 969}%
    \special{pa 2328 906}%
    \special{pa 2359 906}%
    \special{fp}%
    \graphtemp=.5ex\advance\graphtemp by 1.156in
    \rlap{\kern 2.406in\lower\graphtemp\hbox to 0pt{$r_2$\hss}}%
    \graphtemp=.5ex\advance\graphtemp by 1.531in
    \rlap{\kern 2.406in\lower\graphtemp\hbox to 0pt{$s_2$\hss}}%
    \graphtemp=.5ex\advance\graphtemp by 1.906in
    \rlap{\kern 2.438in\lower\graphtemp\hbox to 0pt{$T_2$\hss}}%
    \special{pa 2375 1844}%
    \special{pa 2406 1844}%
    \special{fp}%
    \special{pa 2406 1844}%
    \special{pa 2406 1969}%
    \special{fp}%
    \special{pa 2406 1969}%
    \special{pa 2375 1969}%
    \special{fp}%
    \graphtemp=\baselineskip\multiply\graphtemp by -1\divide\graphtemp by 2
    \advance\graphtemp by .5ex\advance\graphtemp by 0.375in
    \rlap{\kern 1.969in\lower\graphtemp\hbox to 0pt{\hss $Q_3$\hss}}%
    \special{pa 1969 406}%
    \special{pa 1969 594}%
    \special{fp}%
    \special{sh 1.000}%
    \special{pa 1984 531}%
    \special{pa 1969 594}%
    \special{pa 1953 531}%
    \special{pa 1984 531}%
    \special{fp}%
    \graphtemp=.5ex\advance\graphtemp by 0.781in
    \rlap{\kern 2.031in\lower\graphtemp\hbox to 0pt{$r_3$\hss}}%
    \graphtemp=.5ex\advance\graphtemp by 1.156in
    \rlap{\kern 2.031in\lower\graphtemp\hbox to 0pt{$s_3$\hss}}%
    \graphtemp=.5ex\advance\graphtemp by 1.719in
    \rlap{\kern 2.063in\lower\graphtemp\hbox to 0pt{$T_3$\hss}}%
    \special{pa 2000 1469}%
    \special{pa 2031 1469}%
    \special{fp}%
    \special{pa 2031 1469}%
    \special{pa 2031 1969}%
    \special{fp}%
    \special{pa 2031 1969}%
    \special{pa 2000 1969}%
    \special{fp}%
    \graphtemp=\baselineskip\multiply\graphtemp by -1\divide\graphtemp by 2
    \advance\graphtemp by .5ex\advance\graphtemp by 2.656in
    \rlap{\kern 4.031in\lower\graphtemp\hbox to 0pt{\hss $P_4$\hss}}%
    \special{sh 1.000}%
    \special{ar 3469 1906 31 31 0 6.28319}%
    \special{sh 1.000}%
    \special{ar 3469 1531 31 31 0 6.28319}%
    \special{sh 1.000}%
    \special{ar 3469 1156 31 31 0 6.28319}%
    \special{sh 1.000}%
    \special{ar 3469 781 31 31 0 6.28319}%
    \special{sh 1.000}%
    \special{ar 3469 406 31 31 0 6.28319}%
    \special{sh 1.000}%
    \special{ar 3844 1906 31 31 0 6.28319}%
    \special{sh 1.000}%
    \special{ar 3844 1531 31 31 0 6.28319}%
    \special{sh 1.000}%
    \special{ar 3844 1156 31 31 0 6.28319}%
    \special{sh 1.000}%
    \special{ar 3844 781 31 31 0 6.28319}%
    \special{sh 1.000}%
    \special{ar 4219 1906 31 31 0 6.28319}%
    \special{sh 1.000}%
    \special{ar 4219 1531 31 31 0 6.28319}%
    \special{sh 1.000}%
    \special{ar 4219 1156 31 31 0 6.28319}%
    \special{sh 1.000}%
    \special{ar 4594 2281 31 31 0 6.28319}%
    \special{sh 1.000}%
    \special{ar 4594 1906 31 31 0 6.28319}%
    \special{sh 1.000}%
    \special{ar 4594 1531 31 31 0 6.28319}%
    \special{pa 3469 1906}%
    \special{pa 3469 406}%
    \special{fp}%
    \special{pa 3844 1906}%
    \special{pa 3844 781}%
    \special{fp}%
    \special{pa 4219 1906}%
    \special{pa 4219 1156}%
    \special{fp}%
    \special{pa 4594 2281}%
    \special{pa 4594 1531}%
    \special{fp}%
    \special{pa 3469 781}%
    \special{pa 4594 1906}%
    \special{fp}%
    \graphtemp=\baselineskip\multiply\graphtemp by -1\divide\graphtemp by 2
    \advance\graphtemp by .5ex\advance\graphtemp by 1.125in
    \rlap{\kern 4.594in\lower\graphtemp\hbox to 0pt{\hss $Q_1$\hss}}%
    \special{pa 4594 1156}%
    \special{pa 4594 1344}%
    \special{fp}%
    \special{sh 1.000}%
    \special{pa 4609 1281}%
    \special{pa 4594 1344}%
    \special{pa 4578 1281}%
    \special{pa 4609 1281}%
    \special{fp}%
    \graphtemp=.5ex\advance\graphtemp by 1.531in
    \rlap{\kern 4.656in\lower\graphtemp\hbox to 0pt{$r_1$\hss}}%
    \graphtemp=.5ex\advance\graphtemp by 1.906in
    \rlap{\kern 4.656in\lower\graphtemp\hbox to 0pt{$s_1$\hss}}%
    \graphtemp=.5ex\advance\graphtemp by 2.281in
    \rlap{\kern 4.656in\lower\graphtemp\hbox to 0pt{$u$\hss}}%
    \graphtemp=.5ex\advance\graphtemp by 2.281in
    \rlap{\kern 4.500in\lower\graphtemp\hbox to 0pt{\hss $T_1$}}%
    \special{pa 4563 2219}%
    \special{pa 4531 2219}%
    \special{fp}%
    \special{pa 4531 2219}%
    \special{pa 4531 2344}%
    \special{fp}%
    \special{pa 4531 2344}%
    \special{pa 4563 2344}%
    \special{fp}%
    \graphtemp=\baselineskip\multiply\graphtemp by -1\divide\graphtemp by 2
    \advance\graphtemp by .5ex\advance\graphtemp by 0.750in
    \rlap{\kern 4.219in\lower\graphtemp\hbox to 0pt{\hss $Q_2$\hss}}%
    \special{pa 4219 781}%
    \special{pa 4219 969}%
    \special{fp}%
    \special{sh 1.000}%
    \special{pa 4234 906}%
    \special{pa 4219 969}%
    \special{pa 4203 906}%
    \special{pa 4234 906}%
    \special{fp}%
    \graphtemp=.5ex\advance\graphtemp by 1.156in
    \rlap{\kern 4.281in\lower\graphtemp\hbox to 0pt{$r_2$\hss}}%
    \graphtemp=.5ex\advance\graphtemp by 1.531in
    \rlap{\kern 4.281in\lower\graphtemp\hbox to 0pt{$s_2$\hss}}%
    \graphtemp=.5ex\advance\graphtemp by 1.906in
    \rlap{\kern 4.313in\lower\graphtemp\hbox to 0pt{$T_2$\hss}}%
    \special{pa 4250 1844}%
    \special{pa 4281 1844}%
    \special{fp}%
    \special{pa 4281 1844}%
    \special{pa 4281 1969}%
    \special{fp}%
    \special{pa 4281 1969}%
    \special{pa 4250 1969}%
    \special{fp}%
    \graphtemp=\baselineskip\multiply\graphtemp by -1\divide\graphtemp by 2
    \advance\graphtemp by .5ex\advance\graphtemp by 0.375in
    \rlap{\kern 3.844in\lower\graphtemp\hbox to 0pt{\hss $Q_3$\hss}}%
    \special{pa 3844 406}%
    \special{pa 3844 594}%
    \special{fp}%
    \special{sh 1.000}%
    \special{pa 3859 531}%
    \special{pa 3844 594}%
    \special{pa 3828 531}%
    \special{pa 3859 531}%
    \special{fp}%
    \graphtemp=.5ex\advance\graphtemp by 0.781in
    \rlap{\kern 3.906in\lower\graphtemp\hbox to 0pt{$r_3$\hss}}%
    \graphtemp=.5ex\advance\graphtemp by 1.156in
    \rlap{\kern 3.906in\lower\graphtemp\hbox to 0pt{$s_3$\hss}}%
    \graphtemp=.5ex\advance\graphtemp by 1.719in
    \rlap{\kern 3.938in\lower\graphtemp\hbox to 0pt{$T_3$\hss}}%
    \special{pa 3875 1469}%
    \special{pa 3906 1469}%
    \special{fp}%
    \special{pa 3906 1469}%
    \special{pa 3906 1969}%
    \special{fp}%
    \special{pa 3906 1969}%
    \special{pa 3875 1969}%
    \special{fp}%
    \graphtemp=\baselineskip\multiply\graphtemp by -1\divide\graphtemp by 2
    \advance\graphtemp by .5ex\advance\graphtemp by 0.000in
    \rlap{\kern 3.469in\lower\graphtemp\hbox to 0pt{\hss $Q_4$\hss}}%
    \special{pa 3469 31}%
    \special{pa 3469 219}%
    \special{fp}%
    \special{sh 1.000}%
    \special{pa 3484 156}%
    \special{pa 3469 219}%
    \special{pa 3453 156}%
    \special{pa 3484 156}%
    \special{fp}%
    \graphtemp=.5ex\advance\graphtemp by 0.406in
    \rlap{\kern 3.531in\lower\graphtemp\hbox to 0pt{$r_4$\hss}}%
    \graphtemp=.5ex\advance\graphtemp by 0.781in
    \rlap{\kern 3.531in\lower\graphtemp\hbox to 0pt{$s_4$\hss}}%
    \graphtemp=.5ex\advance\graphtemp by 1.531in
    \rlap{\kern 3.563in\lower\graphtemp\hbox to 0pt{$T_4$\hss}}%
    \special{pa 3500 1094}%
    \special{pa 3531 1094}%
    \special{fp}%
    \special{pa 3531 1094}%
    \special{pa 3531 1969}%
    \special{fp}%
    \special{pa 3531 1969}%
    \special{pa 3500 1969}%
    \special{fp}%
    \hbox{\vrule depth2.656in width0pt height 0pt}%
    \kern 4.656in
  }%
}%
\ht\graph=\dp\graph
\dp\graph=0pt

\[\box\graph\]
\caption{The first four $P_j$'s.}
\label{F:labeled}
\end{figure}

\begin{definition}[posets $P_j$] \label{D:mine}
We define a sequence of posets inductively.
Let $P_1$ be a chain $Q_1$ of three elements, $u<s_1<r_1$,
and let $T_1=\left\{u\right\}$
(see Figure~\ref{F:labeled}).
For $j>1$, suppose $P_{j-1}$ is defined and contains the
element $s_{j-1}$.
We define $P_j$ to be the disjoint union of $P_{j-1}$ with a
chain $Q_j$ of $j+1$ elements, plus one new cover relation.
Let $r_j$ be the maximal element of $Q_j$, and let $s_j$
be the next-greatest element.
The set of all remaining elements in $Q_j$ will be called $T_j$.
We add the cover relation $s_{j-1}\lcover s_j$ and extend by
transitivity.\ggcenddef\end{definition}

Figure~\ref{F:labeled} shows $P_1$ through $P_4$.

Next we list some properties of the $P_j$'s.
A ranked poset $P$ has the \emph{strong Sperner property} if, for each
positive integer $k$, the union of the $k$ largest ranks of $P$ is a
maximum-sized $k$-family.
A poset has order dimension at most $2$ if and only if the complement
of its comparability graph is also
a comparability graph~\cite[Thm.~3.61]{DuMi41}
(see also~\cite[p.~62]{TroW92}).

\begin{lemma} \label{L:mineprops}
For the poset $P_j$ of Definition~\ref{D:mine}
the following all hold.
\begin{enumerate}
\item $P_j$ is a ranked poset with width $j$, height $j+2$,
  and cardinality $\binom{j+2}{2}$,
  \label{I:mineprops-ranked}
\item $P_j$ has the strong Sperner property, \label{I:mineprops-stsp}
\item $\Delta d_1(P_j)=j$, $\Delta d_{j+2}(P_j)=1$, and, for
$k=2,\dotsc,j+1$, $\Delta d_k(P_j)=j+2-k$, and
  \label{I:mineprops-delta}
\item $P_j$ has dimension at most $2$. \label{I:mineprops-dim2}
\end{enumerate}
\end{lemma}

\begin{figure}[htbp]
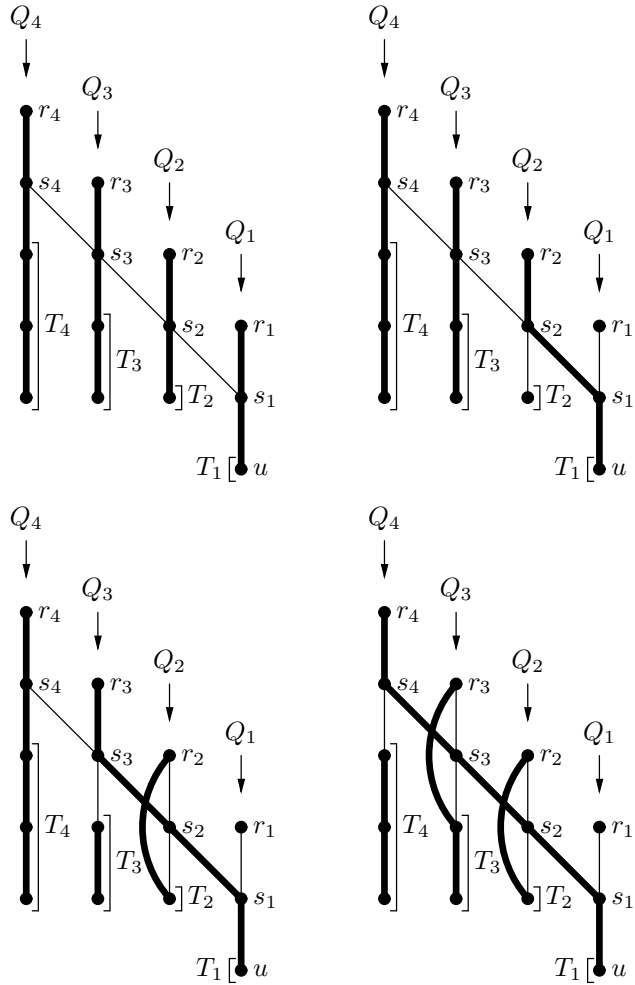

\expandafter\ifx\csname graph\endcsname\relax \csname newbox\endcsname\graph\fi
\expandafter\ifx\csname graphtemp\endcsname\relax \csname newdimen\endcsname\graphtemp\fi
\setbox\graph=\vtop{\vskip 0pt\hbox{%
\fontsize{10}{12}\selectfont 
    \special{pn 8}%
    \special{sh 1.000}%
    \special{ar 31 1906 31 31 0 6.28319}%
    \special{sh 1.000}%
    \special{ar 31 1531 31 31 0 6.28319}%
    \special{sh 1.000}%
    \special{ar 31 1156 31 31 0 6.28319}%
    \special{sh 1.000}%
    \special{ar 31 781 31 31 0 6.28319}%
    \special{sh 1.000}%
    \special{ar 31 406 31 31 0 6.28319}%
    \special{sh 1.000}%
    \special{ar 406 1906 31 31 0 6.28319}%
    \special{sh 1.000}%
    \special{ar 406 1531 31 31 0 6.28319}%
    \special{sh 1.000}%
    \special{ar 406 1156 31 31 0 6.28319}%
    \special{sh 1.000}%
    \special{ar 406 781 31 31 0 6.28319}%
    \special{sh 1.000}%
    \special{ar 781 1906 31 31 0 6.28319}%
    \special{sh 1.000}%
    \special{ar 781 1531 31 31 0 6.28319}%
    \special{sh 1.000}%
    \special{ar 781 1156 31 31 0 6.28319}%
    \special{sh 1.000}%
    \special{ar 1156 2281 31 31 0 6.28319}%
    \special{sh 1.000}%
    \special{ar 1156 1906 31 31 0 6.28319}%
    \special{sh 1.000}%
    \special{ar 1156 1531 31 31 0 6.28319}%
    \special{pa 31 1906}%
    \special{pa 31 406}%
    \special{fp}%
    \special{pa 406 1906}%
    \special{pa 406 781}%
    \special{fp}%
    \special{pa 781 1906}%
    \special{pa 781 1156}%
    \special{fp}%
    \special{pa 1156 2281}%
    \special{pa 1156 1531}%
    \special{fp}%
    \special{pa 31 781}%
    \special{pa 1156 1906}%
    \special{fp}%
    \graphtemp=\baselineskip\multiply\graphtemp by -1\divide\graphtemp by 2
    \advance\graphtemp by .5ex\advance\graphtemp by 1.125in
    \rlap{\kern 1.156in\lower\graphtemp\hbox to 0pt{\hss $Q_1$\hss}}%
    \special{pa 1156 1156}%
    \special{pa 1156 1344}%
    \special{fp}%
    \special{sh 1.000}%
    \special{pa 1172 1281}%
    \special{pa 1156 1344}%
    \special{pa 1141 1281}%
    \special{pa 1172 1281}%
    \special{fp}%
    \graphtemp=.5ex\advance\graphtemp by 1.531in
    \rlap{\kern 1.219in\lower\graphtemp\hbox to 0pt{$r_1$\hss}}%
    \graphtemp=.5ex\advance\graphtemp by 1.906in
    \rlap{\kern 1.219in\lower\graphtemp\hbox to 0pt{$s_1$\hss}}%
    \graphtemp=.5ex\advance\graphtemp by 2.281in
    \rlap{\kern 1.219in\lower\graphtemp\hbox to 0pt{$u$\hss}}%
    \graphtemp=.5ex\advance\graphtemp by 2.281in
    \rlap{\kern 1.063in\lower\graphtemp\hbox to 0pt{\hss $T_1$}}%
    \special{pa 1125 2219}%
    \special{pa 1094 2219}%
    \special{fp}%
    \special{pa 1094 2219}%
    \special{pa 1094 2344}%
    \special{fp}%
    \special{pa 1094 2344}%
    \special{pa 1125 2344}%
    \special{fp}%
    \graphtemp=\baselineskip\multiply\graphtemp by -1\divide\graphtemp by 2
    \advance\graphtemp by .5ex\advance\graphtemp by 0.750in
    \rlap{\kern 0.781in\lower\graphtemp\hbox to 0pt{\hss $Q_2$\hss}}%
    \special{pa 781 781}%
    \special{pa 781 969}%
    \special{fp}%
    \special{sh 1.000}%
    \special{pa 797 906}%
    \special{pa 781 969}%
    \special{pa 766 906}%
    \special{pa 797 906}%
    \special{fp}%
    \graphtemp=.5ex\advance\graphtemp by 1.156in
    \rlap{\kern 0.844in\lower\graphtemp\hbox to 0pt{$r_2$\hss}}%
    \graphtemp=.5ex\advance\graphtemp by 1.531in
    \rlap{\kern 0.844in\lower\graphtemp\hbox to 0pt{$s_2$\hss}}%
    \graphtemp=.5ex\advance\graphtemp by 1.906in
    \rlap{\kern 0.875in\lower\graphtemp\hbox to 0pt{$T_2$\hss}}%
    \special{pa 813 1844}%
    \special{pa 844 1844}%
    \special{fp}%
    \special{pa 844 1844}%
    \special{pa 844 1969}%
    \special{fp}%
    \special{pa 844 1969}%
    \special{pa 813 1969}%
    \special{fp}%
    \graphtemp=\baselineskip\multiply\graphtemp by -1\divide\graphtemp by 2
    \advance\graphtemp by .5ex\advance\graphtemp by 0.375in
    \rlap{\kern 0.406in\lower\graphtemp\hbox to 0pt{\hss $Q_3$\hss}}%
    \special{pa 406 406}%
    \special{pa 406 594}%
    \special{fp}%
    \special{sh 1.000}%
    \special{pa 422 531}%
    \special{pa 406 594}%
    \special{pa 391 531}%
    \special{pa 422 531}%
    \special{fp}%
    \graphtemp=.5ex\advance\graphtemp by 0.781in
    \rlap{\kern 0.469in\lower\graphtemp\hbox to 0pt{$r_3$\hss}}%
    \graphtemp=.5ex\advance\graphtemp by 1.156in
    \rlap{\kern 0.469in\lower\graphtemp\hbox to 0pt{$s_3$\hss}}%
    \graphtemp=.5ex\advance\graphtemp by 1.719in
    \rlap{\kern 0.500in\lower\graphtemp\hbox to 0pt{$T_3$\hss}}%
    \special{pa 438 1469}%
    \special{pa 469 1469}%
    \special{fp}%
    \special{pa 469 1469}%
    \special{pa 469 1969}%
    \special{fp}%
    \special{pa 469 1969}%
    \special{pa 438 1969}%
    \special{fp}%
    \graphtemp=\baselineskip\multiply\graphtemp by -1\divide\graphtemp by 2
    \advance\graphtemp by .5ex\advance\graphtemp by 0.000in
    \rlap{\kern 0.031in\lower\graphtemp\hbox to 0pt{\hss $Q_4$\hss}}%
    \special{pa 31 31}%
    \special{pa 31 219}%
    \special{fp}%
    \special{sh 1.000}%
    \special{pa 47 156}%
    \special{pa 31 219}%
    \special{pa 16 156}%
    \special{pa 47 156}%
    \special{fp}%
    \graphtemp=.5ex\advance\graphtemp by 0.406in
    \rlap{\kern 0.094in\lower\graphtemp\hbox to 0pt{$r_4$\hss}}%
    \graphtemp=.5ex\advance\graphtemp by 0.781in
    \rlap{\kern 0.094in\lower\graphtemp\hbox to 0pt{$s_4$\hss}}%
    \graphtemp=.5ex\advance\graphtemp by 1.531in
    \rlap{\kern 0.125in\lower\graphtemp\hbox to 0pt{$T_4$\hss}}%
    \special{pa 63 1094}%
    \special{pa 94 1094}%
    \special{fp}%
    \special{pa 94 1094}%
    \special{pa 94 1969}%
    \special{fp}%
    \special{pa 94 1969}%
    \special{pa 63 1969}%
    \special{fp}%
    \special{pn 35}%
    \special{pa 31 1906}%
    \special{pa 31 406}%
    \special{fp}%
    \special{pa 406 1906}%
    \special{pa 406 781}%
    \special{fp}%
    \special{pa 781 1906}%
    \special{pa 781 1156}%
    \special{fp}%
    \special{pa 1156 2281}%
    \special{pa 1156 1531}%
    \special{fp}%
    \special{pn 8}%
    \special{sh 1.000}%
    \special{ar 1906 1906 31 31 0 6.28319}%
    \special{sh 1.000}%
    \special{ar 1906 1531 31 31 0 6.28319}%
    \special{sh 1.000}%
    \special{ar 1906 1156 31 31 0 6.28319}%
    \special{sh 1.000}%
    \special{ar 1906 781 31 31 0 6.28319}%
    \special{sh 1.000}%
    \special{ar 1906 406 31 31 0 6.28319}%
    \special{sh 1.000}%
    \special{ar 2281 1906 31 31 0 6.28319}%
    \special{sh 1.000}%
    \special{ar 2281 1531 31 31 0 6.28319}%
    \special{sh 1.000}%
    \special{ar 2281 1156 31 31 0 6.28319}%
    \special{sh 1.000}%
    \special{ar 2281 781 31 31 0 6.28319}%
    \special{sh 1.000}%
    \special{ar 2656 1906 31 31 0 6.28319}%
    \special{sh 1.000}%
    \special{ar 2656 1531 31 31 0 6.28319}%
    \special{sh 1.000}%
    \special{ar 2656 1156 31 31 0 6.28319}%
    \special{sh 1.000}%
    \special{ar 3031 2281 31 31 0 6.28319}%
    \special{sh 1.000}%
    \special{ar 3031 1906 31 31 0 6.28319}%
    \special{sh 1.000}%
    \special{ar 3031 1531 31 31 0 6.28319}%
    \special{pa 1906 1906}%
    \special{pa 1906 406}%
    \special{fp}%
    \special{pa 2281 1906}%
    \special{pa 2281 781}%
    \special{fp}%
    \special{pa 2656 1906}%
    \special{pa 2656 1156}%
    \special{fp}%
    \special{pa 3031 2281}%
    \special{pa 3031 1531}%
    \special{fp}%
    \special{pa 1906 781}%
    \special{pa 3031 1906}%
    \special{fp}%
    \graphtemp=\baselineskip\multiply\graphtemp by -1\divide\graphtemp by 2
    \advance\graphtemp by .5ex\advance\graphtemp by 1.125in
    \rlap{\kern 3.031in\lower\graphtemp\hbox to 0pt{\hss $Q_1$\hss}}%
    \special{pa 3031 1156}%
    \special{pa 3031 1344}%
    \special{fp}%
    \special{sh 1.000}%
    \special{pa 3047 1281}%
    \special{pa 3031 1344}%
    \special{pa 3016 1281}%
    \special{pa 3047 1281}%
    \special{fp}%
    \graphtemp=.5ex\advance\graphtemp by 1.531in
    \rlap{\kern 3.094in\lower\graphtemp\hbox to 0pt{$r_1$\hss}}%
    \graphtemp=.5ex\advance\graphtemp by 1.906in
    \rlap{\kern 3.094in\lower\graphtemp\hbox to 0pt{$s_1$\hss}}%
    \graphtemp=.5ex\advance\graphtemp by 2.281in
    \rlap{\kern 3.094in\lower\graphtemp\hbox to 0pt{$u$\hss}}%
    \graphtemp=.5ex\advance\graphtemp by 2.281in
    \rlap{\kern 2.938in\lower\graphtemp\hbox to 0pt{\hss $T_1$}}%
    \special{pa 3000 2219}%
    \special{pa 2969 2219}%
    \special{fp}%
    \special{pa 2969 2219}%
    \special{pa 2969 2344}%
    \special{fp}%
    \special{pa 2969 2344}%
    \special{pa 3000 2344}%
    \special{fp}%
    \graphtemp=\baselineskip\multiply\graphtemp by -1\divide\graphtemp by 2
    \advance\graphtemp by .5ex\advance\graphtemp by 0.750in
    \rlap{\kern 2.656in\lower\graphtemp\hbox to 0pt{\hss $Q_2$\hss}}%
    \special{pa 2656 781}%
    \special{pa 2656 969}%
    \special{fp}%
    \special{sh 1.000}%
    \special{pa 2672 906}%
    \special{pa 2656 969}%
    \special{pa 2641 906}%
    \special{pa 2672 906}%
    \special{fp}%
    \graphtemp=.5ex\advance\graphtemp by 1.156in
    \rlap{\kern 2.719in\lower\graphtemp\hbox to 0pt{$r_2$\hss}}%
    \graphtemp=.5ex\advance\graphtemp by 1.531in
    \rlap{\kern 2.719in\lower\graphtemp\hbox to 0pt{$s_2$\hss}}%
    \graphtemp=.5ex\advance\graphtemp by 1.906in
    \rlap{\kern 2.750in\lower\graphtemp\hbox to 0pt{$T_2$\hss}}%
    \special{pa 2688 1844}%
    \special{pa 2719 1844}%
    \special{fp}%
    \special{pa 2719 1844}%
    \special{pa 2719 1969}%
    \special{fp}%
    \special{pa 2719 1969}%
    \special{pa 2688 1969}%
    \special{fp}%
    \graphtemp=\baselineskip\multiply\graphtemp by -1\divide\graphtemp by 2
    \advance\graphtemp by .5ex\advance\graphtemp by 0.375in
    \rlap{\kern 2.281in\lower\graphtemp\hbox to 0pt{\hss $Q_3$\hss}}%
    \special{pa 2281 406}%
    \special{pa 2281 594}%
    \special{fp}%
    \special{sh 1.000}%
    \special{pa 2297 531}%
    \special{pa 2281 594}%
    \special{pa 2266 531}%
    \special{pa 2297 531}%
    \special{fp}%
    \graphtemp=.5ex\advance\graphtemp by 0.781in
    \rlap{\kern 2.344in\lower\graphtemp\hbox to 0pt{$r_3$\hss}}%
    \graphtemp=.5ex\advance\graphtemp by 1.156in
    \rlap{\kern 2.344in\lower\graphtemp\hbox to 0pt{$s_3$\hss}}%
    \graphtemp=.5ex\advance\graphtemp by 1.719in
    \rlap{\kern 2.375in\lower\graphtemp\hbox to 0pt{$T_3$\hss}}%
    \special{pa 2313 1469}%
    \special{pa 2344 1469}%
    \special{fp}%
    \special{pa 2344 1469}%
    \special{pa 2344 1969}%
    \special{fp}%
    \special{pa 2344 1969}%
    \special{pa 2313 1969}%
    \special{fp}%
    \graphtemp=\baselineskip\multiply\graphtemp by -1\divide\graphtemp by 2
    \advance\graphtemp by .5ex\advance\graphtemp by 0.000in
    \rlap{\kern 1.906in\lower\graphtemp\hbox to 0pt{\hss $Q_4$\hss}}%
    \special{pa 1906 31}%
    \special{pa 1906 219}%
    \special{fp}%
    \special{sh 1.000}%
    \special{pa 1922 156}%
    \special{pa 1906 219}%
    \special{pa 1891 156}%
    \special{pa 1922 156}%
    \special{fp}%
    \graphtemp=.5ex\advance\graphtemp by 0.406in
    \rlap{\kern 1.969in\lower\graphtemp\hbox to 0pt{$r_4$\hss}}%
    \graphtemp=.5ex\advance\graphtemp by 0.781in
    \rlap{\kern 1.969in\lower\graphtemp\hbox to 0pt{$s_4$\hss}}%
    \graphtemp=.5ex\advance\graphtemp by 1.531in
    \rlap{\kern 2.000in\lower\graphtemp\hbox to 0pt{$T_4$\hss}}%
    \special{pa 1938 1094}%
    \special{pa 1969 1094}%
    \special{fp}%
    \special{pa 1969 1094}%
    \special{pa 1969 1969}%
    \special{fp}%
    \special{pa 1969 1969}%
    \special{pa 1938 1969}%
    \special{fp}%
    \special{pn 35}%
    \special{pa 1906 1906}%
    \special{pa 1906 406}%
    \special{fp}%
    \special{pa 2281 1906}%
    \special{pa 2281 781}%
    \special{fp}%
    \special{pa 2656 1156}%
    \special{pa 2656 1531}%
    \special{fp}%
    \special{pa 2656 1531}%
    \special{pa 3031 1906}%
    \special{fp}%
    \special{pa 3031 1906}%
    \special{pa 3031 2281}%
    \special{fp}%
    \special{pn 8}%
    \special{sh 1.000}%
    \special{ar 31 4531 31 31 0 6.28319}%
    \special{sh 1.000}%
    \special{ar 31 4156 31 31 0 6.28319}%
    \special{sh 1.000}%
    \special{ar 31 3781 31 31 0 6.28319}%
    \special{sh 1.000}%
    \special{ar 31 3406 31 31 0 6.28319}%
    \special{sh 1.000}%
    \special{ar 31 3031 31 31 0 6.28319}%
    \special{sh 1.000}%
    \special{ar 406 4531 31 31 0 6.28319}%
    \special{sh 1.000}%
    \special{ar 406 4156 31 31 0 6.28319}%
    \special{sh 1.000}%
    \special{ar 406 3781 31 31 0 6.28319}%
    \special{sh 1.000}%
    \special{ar 406 3406 31 31 0 6.28319}%
    \special{sh 1.000}%
    \special{ar 781 4531 31 31 0 6.28319}%
    \special{sh 1.000}%
    \special{ar 781 4156 31 31 0 6.28319}%
    \special{sh 1.000}%
    \special{ar 781 3781 31 31 0 6.28319}%
    \special{sh 1.000}%
    \special{ar 1156 4906 31 31 0 6.28319}%
    \special{sh 1.000}%
    \special{ar 1156 4531 31 31 0 6.28319}%
    \special{sh 1.000}%
    \special{ar 1156 4156 31 31 0 6.28319}%
    \special{pa 31 4531}%
    \special{pa 31 3031}%
    \special{fp}%
    \special{pa 406 4531}%
    \special{pa 406 3406}%
    \special{fp}%
    \special{pa 781 4531}%
    \special{pa 781 3781}%
    \special{fp}%
    \special{pa 1156 4906}%
    \special{pa 1156 4156}%
    \special{fp}%
    \special{pa 31 3406}%
    \special{pa 1156 4531}%
    \special{fp}%
    \graphtemp=\baselineskip\multiply\graphtemp by -1\divide\graphtemp by 2
    \advance\graphtemp by .5ex\advance\graphtemp by 3.750in
    \rlap{\kern 1.156in\lower\graphtemp\hbox to 0pt{\hss $Q_1$\hss}}%
    \special{pa 1156 3781}%
    \special{pa 1156 3969}%
    \special{fp}%
    \special{sh 1.000}%
    \special{pa 1172 3906}%
    \special{pa 1156 3969}%
    \special{pa 1141 3906}%
    \special{pa 1172 3906}%
    \special{fp}%
    \graphtemp=.5ex\advance\graphtemp by 4.156in
    \rlap{\kern 1.219in\lower\graphtemp\hbox to 0pt{$r_1$\hss}}%
    \graphtemp=.5ex\advance\graphtemp by 4.531in
    \rlap{\kern 1.219in\lower\graphtemp\hbox to 0pt{$s_1$\hss}}%
    \graphtemp=.5ex\advance\graphtemp by 4.906in
    \rlap{\kern 1.219in\lower\graphtemp\hbox to 0pt{$u$\hss}}%
    \graphtemp=.5ex\advance\graphtemp by 4.906in
    \rlap{\kern 1.063in\lower\graphtemp\hbox to 0pt{\hss $T_1$}}%
    \special{pa 1125 4844}%
    \special{pa 1094 4844}%
    \special{fp}%
    \special{pa 1094 4844}%
    \special{pa 1094 4969}%
    \special{fp}%
    \special{pa 1094 4969}%
    \special{pa 1125 4969}%
    \special{fp}%
    \graphtemp=\baselineskip\multiply\graphtemp by -1\divide\graphtemp by 2
    \advance\graphtemp by .5ex\advance\graphtemp by 3.375in
    \rlap{\kern 0.781in\lower\graphtemp\hbox to 0pt{\hss $Q_2$\hss}}%
    \special{pa 781 3406}%
    \special{pa 781 3594}%
    \special{fp}%
    \special{sh 1.000}%
    \special{pa 797 3531}%
    \special{pa 781 3594}%
    \special{pa 766 3531}%
    \special{pa 797 3531}%
    \special{fp}%
    \graphtemp=.5ex\advance\graphtemp by 3.781in
    \rlap{\kern 0.844in\lower\graphtemp\hbox to 0pt{$r_2$\hss}}%
    \graphtemp=.5ex\advance\graphtemp by 4.156in
    \rlap{\kern 0.844in\lower\graphtemp\hbox to 0pt{$s_2$\hss}}%
    \graphtemp=.5ex\advance\graphtemp by 4.531in
    \rlap{\kern 0.875in\lower\graphtemp\hbox to 0pt{$T_2$\hss}}%
    \special{pa 813 4469}%
    \special{pa 844 4469}%
    \special{fp}%
    \special{pa 844 4469}%
    \special{pa 844 4594}%
    \special{fp}%
    \special{pa 844 4594}%
    \special{pa 813 4594}%
    \special{fp}%
    \graphtemp=\baselineskip\multiply\graphtemp by -1\divide\graphtemp by 2
    \advance\graphtemp by .5ex\advance\graphtemp by 3.000in
    \rlap{\kern 0.406in\lower\graphtemp\hbox to 0pt{\hss $Q_3$\hss}}%
    \special{pa 406 3031}%
    \special{pa 406 3219}%
    \special{fp}%
    \special{sh 1.000}%
    \special{pa 422 3156}%
    \special{pa 406 3219}%
    \special{pa 391 3156}%
    \special{pa 422 3156}%
    \special{fp}%
    \graphtemp=.5ex\advance\graphtemp by 3.406in
    \rlap{\kern 0.469in\lower\graphtemp\hbox to 0pt{$r_3$\hss}}%
    \graphtemp=.5ex\advance\graphtemp by 3.781in
    \rlap{\kern 0.469in\lower\graphtemp\hbox to 0pt{$s_3$\hss}}%
    \graphtemp=.5ex\advance\graphtemp by 4.344in
    \rlap{\kern 0.500in\lower\graphtemp\hbox to 0pt{$T_3$\hss}}%
    \special{pa 438 4094}%
    \special{pa 469 4094}%
    \special{fp}%
    \special{pa 469 4094}%
    \special{pa 469 4594}%
    \special{fp}%
    \special{pa 469 4594}%
    \special{pa 438 4594}%
    \special{fp}%
    \graphtemp=\baselineskip\multiply\graphtemp by -1\divide\graphtemp by 2
    \advance\graphtemp by .5ex\advance\graphtemp by 2.625in
    \rlap{\kern 0.031in\lower\graphtemp\hbox to 0pt{\hss $Q_4$\hss}}%
    \special{pa 31 2656}%
    \special{pa 31 2844}%
    \special{fp}%
    \special{sh 1.000}%
    \special{pa 47 2781}%
    \special{pa 31 2844}%
    \special{pa 16 2781}%
    \special{pa 47 2781}%
    \special{fp}%
    \graphtemp=.5ex\advance\graphtemp by 3.031in
    \rlap{\kern 0.094in\lower\graphtemp\hbox to 0pt{$r_4$\hss}}%
    \graphtemp=.5ex\advance\graphtemp by 3.406in
    \rlap{\kern 0.094in\lower\graphtemp\hbox to 0pt{$s_4$\hss}}%
    \graphtemp=.5ex\advance\graphtemp by 4.156in
    \rlap{\kern 0.125in\lower\graphtemp\hbox to 0pt{$T_4$\hss}}%
    \special{pa 63 3719}%
    \special{pa 94 3719}%
    \special{fp}%
    \special{pa 94 3719}%
    \special{pa 94 4594}%
    \special{fp}%
    \special{pa 94 4594}%
    \special{pa 63 4594}%
    \special{fp}%
    \special{pn 35}%
    \special{pa 31 4531}%
    \special{pa 31 3031}%
    \special{fp}%
    \special{pa 406 4531}%
    \special{pa 406 4156}%
    \special{fp}%
    \special{pa 406 3406}%
    \special{pa 406 3781}%
    \special{fp}%
    \special{pa 406 3781}%
    \special{pa 1156 4531}%
    \special{fp}%
    \special{pa 1156 4531}%
    \special{pa 1156 4906}%
    \special{fp}%
    \special{ar 1201 4156 563 563 2.411865 3.871320}%
    \special{pn 8}%
    \special{sh 1.000}%
    \special{ar 1906 4531 31 31 0 6.28319}%
    \special{sh 1.000}%
    \special{ar 1906 4156 31 31 0 6.28319}%
    \special{sh 1.000}%
    \special{ar 1906 3781 31 31 0 6.28319}%
    \special{sh 1.000}%
    \special{ar 1906 3406 31 31 0 6.28319}%
    \special{sh 1.000}%
    \special{ar 1906 3031 31 31 0 6.28319}%
    \special{sh 1.000}%
    \special{ar 2281 4531 31 31 0 6.28319}%
    \special{sh 1.000}%
    \special{ar 2281 4156 31 31 0 6.28319}%
    \special{sh 1.000}%
    \special{ar 2281 3781 31 31 0 6.28319}%
    \special{sh 1.000}%
    \special{ar 2281 3406 31 31 0 6.28319}%
    \special{sh 1.000}%
    \special{ar 2656 4531 31 31 0 6.28319}%
    \special{sh 1.000}%
    \special{ar 2656 4156 31 31 0 6.28319}%
    \special{sh 1.000}%
    \special{ar 2656 3781 31 31 0 6.28319}%
    \special{sh 1.000}%
    \special{ar 3031 4906 31 31 0 6.28319}%
    \special{sh 1.000}%
    \special{ar 3031 4531 31 31 0 6.28319}%
    \special{sh 1.000}%
    \special{ar 3031 4156 31 31 0 6.28319}%
    \special{pa 1906 4531}%
    \special{pa 1906 3031}%
    \special{fp}%
    \special{pa 2281 4531}%
    \special{pa 2281 3406}%
    \special{fp}%
    \special{pa 2656 4531}%
    \special{pa 2656 3781}%
    \special{fp}%
    \special{pa 3031 4906}%
    \special{pa 3031 4156}%
    \special{fp}%
    \special{pa 1906 3406}%
    \special{pa 3031 4531}%
    \special{fp}%
    \graphtemp=\baselineskip\multiply\graphtemp by -1\divide\graphtemp by 2
    \advance\graphtemp by .5ex\advance\graphtemp by 3.750in
    \rlap{\kern 3.031in\lower\graphtemp\hbox to 0pt{\hss $Q_1$\hss}}%
    \special{pa 3031 3781}%
    \special{pa 3031 3969}%
    \special{fp}%
    \special{sh 1.000}%
    \special{pa 3047 3906}%
    \special{pa 3031 3969}%
    \special{pa 3016 3906}%
    \special{pa 3047 3906}%
    \special{fp}%
    \graphtemp=.5ex\advance\graphtemp by 4.156in
    \rlap{\kern 3.094in\lower\graphtemp\hbox to 0pt{$r_1$\hss}}%
    \graphtemp=.5ex\advance\graphtemp by 4.531in
    \rlap{\kern 3.094in\lower\graphtemp\hbox to 0pt{$s_1$\hss}}%
    \graphtemp=.5ex\advance\graphtemp by 4.906in
    \rlap{\kern 3.094in\lower\graphtemp\hbox to 0pt{$u$\hss}}%
    \graphtemp=.5ex\advance\graphtemp by 4.906in
    \rlap{\kern 2.938in\lower\graphtemp\hbox to 0pt{\hss $T_1$}}%
    \special{pa 3000 4844}%
    \special{pa 2969 4844}%
    \special{fp}%
    \special{pa 2969 4844}%
    \special{pa 2969 4969}%
    \special{fp}%
    \special{pa 2969 4969}%
    \special{pa 3000 4969}%
    \special{fp}%
    \graphtemp=\baselineskip\multiply\graphtemp by -1\divide\graphtemp by 2
    \advance\graphtemp by .5ex\advance\graphtemp by 3.375in
    \rlap{\kern 2.656in\lower\graphtemp\hbox to 0pt{\hss $Q_2$\hss}}%
    \special{pa 2656 3406}%
    \special{pa 2656 3594}%
    \special{fp}%
    \special{sh 1.000}%
    \special{pa 2672 3531}%
    \special{pa 2656 3594}%
    \special{pa 2641 3531}%
    \special{pa 2672 3531}%
    \special{fp}%
    \graphtemp=.5ex\advance\graphtemp by 3.781in
    \rlap{\kern 2.719in\lower\graphtemp\hbox to 0pt{$r_2$\hss}}%
    \graphtemp=.5ex\advance\graphtemp by 4.156in
    \rlap{\kern 2.719in\lower\graphtemp\hbox to 0pt{$s_2$\hss}}%
    \graphtemp=.5ex\advance\graphtemp by 4.531in
    \rlap{\kern 2.750in\lower\graphtemp\hbox to 0pt{$T_2$\hss}}%
    \special{pa 2688 4469}%
    \special{pa 2719 4469}%
    \special{fp}%
    \special{pa 2719 4469}%
    \special{pa 2719 4594}%
    \special{fp}%
    \special{pa 2719 4594}%
    \special{pa 2688 4594}%
    \special{fp}%
    \graphtemp=\baselineskip\multiply\graphtemp by -1\divide\graphtemp by 2
    \advance\graphtemp by .5ex\advance\graphtemp by 3.000in
    \rlap{\kern 2.281in\lower\graphtemp\hbox to 0pt{\hss $Q_3$\hss}}%
    \special{pa 2281 3031}%
    \special{pa 2281 3219}%
    \special{fp}%
    \special{sh 1.000}%
    \special{pa 2297 3156}%
    \special{pa 2281 3219}%
    \special{pa 2266 3156}%
    \special{pa 2297 3156}%
    \special{fp}%
    \graphtemp=.5ex\advance\graphtemp by 3.406in
    \rlap{\kern 2.344in\lower\graphtemp\hbox to 0pt{$r_3$\hss}}%
    \graphtemp=.5ex\advance\graphtemp by 3.781in
    \rlap{\kern 2.344in\lower\graphtemp\hbox to 0pt{$s_3$\hss}}%
    \graphtemp=.5ex\advance\graphtemp by 4.344in
    \rlap{\kern 2.375in\lower\graphtemp\hbox to 0pt{$T_3$\hss}}%
    \special{pa 2313 4094}%
    \special{pa 2344 4094}%
    \special{fp}%
    \special{pa 2344 4094}%
    \special{pa 2344 4594}%
    \special{fp}%
    \special{pa 2344 4594}%
    \special{pa 2313 4594}%
    \special{fp}%
    \graphtemp=\baselineskip\multiply\graphtemp by -1\divide\graphtemp by 2
    \advance\graphtemp by .5ex\advance\graphtemp by 2.625in
    \rlap{\kern 1.906in\lower\graphtemp\hbox to 0pt{\hss $Q_4$\hss}}%
    \special{pa 1906 2656}%
    \special{pa 1906 2844}%
    \special{fp}%
    \special{sh 1.000}%
    \special{pa 1922 2781}%
    \special{pa 1906 2844}%
    \special{pa 1891 2781}%
    \special{pa 1922 2781}%
    \special{fp}%
    \graphtemp=.5ex\advance\graphtemp by 3.031in
    \rlap{\kern 1.969in\lower\graphtemp\hbox to 0pt{$r_4$\hss}}%
    \graphtemp=.5ex\advance\graphtemp by 3.406in
    \rlap{\kern 1.969in\lower\graphtemp\hbox to 0pt{$s_4$\hss}}%
    \graphtemp=.5ex\advance\graphtemp by 4.156in
    \rlap{\kern 2.000in\lower\graphtemp\hbox to 0pt{$T_4$\hss}}%
    \special{pa 1938 3719}%
    \special{pa 1969 3719}%
    \special{fp}%
    \special{pa 1969 3719}%
    \special{pa 1969 4594}%
    \special{fp}%
    \special{pa 1969 4594}%
    \special{pa 1938 4594}%
    \special{fp}%
    \special{pn 35}%
    \special{pa 1906 4531}%
    \special{pa 1906 3781}%
    \special{fp}%
    \special{pa 2281 4531}%
    \special{pa 2281 4156}%
    \special{fp}%
    \special{pa 1906 3031}%
    \special{pa 1906 3406}%
    \special{fp}%
    \special{pa 1906 3406}%
    \special{pa 3031 4531}%
    \special{fp}%
    \special{pa 3031 4531}%
    \special{pa 3031 4906}%
    \special{fp}%
    \special{ar 2701 3781 563 563 2.411865 3.871320}%
    \special{ar 3076 4156 563 563 2.411865 3.871320}%
    \hbox{\vrule depth4.969in width0pt height 0pt}%
    \kern 3.094in
  }%
}%
\ht\graph=\dp\graph
\dp\graph=0pt

\[\box\graph\]
\caption{$P_4$ shown with chain partitions
$\cC_1$, $\cC_2$, $\cC_3$, and $\cC_4$,
as defined in the proof of Lemma~\ref{L:mineprops}.}
\label{F:part}
\end{figure}

\begin{proof}
\ref{I:mineprops-ranked}
This is immediate from the construction, with $u<s_1<\dotsb<s_j<r_j$ a
maximum chain and $\left\{r_i\right\}$ a maximum antichain.

\ref{I:mineprops-stsp}
Let $k$ be a positive integer.
Let $S$ be the union of the $k$ largest
ranks of $P_j$.
It suffices to exhibit a chain partition
$\cC_k$ of $P_j$ such that $|S|=m_k(\cC_k)$.
Let $\cC_k$ consist of the chain
$C=\left\{u,s_1,\dotsc,s_k,r_k\right\}$ and the chains $Q_i\setmin C$
for $1\le i\le k$
(see Figure~\ref{F:part} for examples).
We claim that $\cC_k$ is the required partition.
Each chain $C\in\cC_k$ either is contained in $S$
or contains one element of each rank of $S$.
In the former case, $C$
contributes $|C|$ to $m_k(\cC_k)$; in the latter case, $C$
contributes $k$.
In both cases, $C$ contributes $\min\left\{k,|C|\right\}$,
and so $m_k(\cC_k)=|S|$.
(In fact, $\cC_k$ is both $k$- and $k+1$-saturated.)

\ref{I:mineprops-delta}
This follows from Statement~\ref{I:mineprops-stsp} and the
sequence of rank sizes of $P_j$.

\ref{I:mineprops-dim2}
We define two linear extensions of $P_j$:
\begin{itemize}
\item
  $T_1,s_1,r_1,\quad T_2,s_2,r_2,\quad\dotsc,\quad T_j,s_j,r_j$,
  \label{I:order-byc}
\item
  $T_j,\dotsc,T_2,T_1,\quad s_1,s_2,\dotsc,s_j,\quad r_j,
  \dotsc,r_2,r_1$,
  \label{I:order-dud}
\end{itemize}
where the elements of $T_i$ are in ascending order.
The intersection of
these two total orders is the partial order on $P_j$, i.e., each
related pair appears in the proper order in both extensions, and each
incomparable pair appears in opposite orders in the two extensions.
Thus, $P_j$ has dimension at most $2$.\ggcendpf\end{proof}

\begin{lemma} \label{L:minepp}
For each positive integer $j$, $P_j$ is a polyunsaturated poset.
\end{lemma}

\begin{proof}
We claim that, if $2\le k\le j$ and $\cC$ is a $k$-saturated
chain partition of
$P_j$, then $u$ lies on the same chain in $\cC$ as either
$r_{k-1}$ or $r_k$.
Also, in every $1$-saturated chain partition, $u$
lies on the same chain as $r_1$, and in every $j+1$-saturated chain
partition, $u$ lies on the same chain as $r_j$.
Since in $\cC$
the element $u$ can be on the same chain with only one $r_i$, this
prevents $\cC$ from being $k$-saturated for nonconsecutive
values of $k$.

Let $2\le k\le j$, and
let $\cC$ be a $k$-saturated chain partition of $P_j$.
The ranks containing $s_1,\dotsc,s_j$, in order, are the $j$ largest
ranks.
Thus there exists a union of $k$ largest ranks of $P_j$ that omits
both $u$ and $r_k$.
By Statement~\ref{I:mineprops-stsp} of Lemma~\ref{L:mineprops}, a
union of $k$ largest ranks of $P_j$ is a maximum $k$-family.
Thus, the chain in $\cC$ that contains $u$ must contain at
least one element from each of those ranks, and the chain containing
$r_k$ must contain at least one element from each of those ranks.
Since $r_{k-1}$ and $s_k$ are the only elements of the $k$th largest
rank that are comparable to $u$, one of the two must be on its chain.
Similarly, $s_k$ must be on the same chain as $r_k$ in $\cC$.
Since either $r_{k-1}$ or $s_k$ must be on the chain containing $u$,
we see that either $r_{k-1}$ or $r_k$ must be on this chain.

The $k=1$ and $k=j+1$ cases of the claim are proven similarly.
Thus, we have established the claim.\ggcendpf\end{proof}

\begin{figure}[htbp]
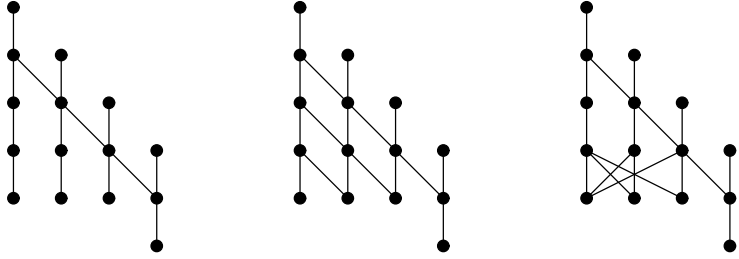

\expandafter\ifx\csname graph\endcsname\relax \csname newbox\endcsname\graph\fi
\expandafter\ifx\csname graphtemp\endcsname\relax \csname newdimen\endcsname\graphtemp\fi
\setbox\graph=\vtop{\vskip 0pt\hbox{%
\fontsize{10}{12}\selectfont 
    \special{pn 8}%
    \special{sh 1.000}%
    \special{ar 31 1031 31 31 0 6.28319}%
    \special{sh 1.000}%
    \special{ar 31 781 31 31 0 6.28319}%
    \special{sh 1.000}%
    \special{ar 31 531 31 31 0 6.28319}%
    \special{sh 1.000}%
    \special{ar 31 281 31 31 0 6.28319}%
    \special{sh 1.000}%
    \special{ar 31 31 31 31 0 6.28319}%
    \special{sh 1.000}%
    \special{ar 281 1031 31 31 0 6.28319}%
    \special{sh 1.000}%
    \special{ar 281 781 31 31 0 6.28319}%
    \special{sh 1.000}%
    \special{ar 281 531 31 31 0 6.28319}%
    \special{sh 1.000}%
    \special{ar 281 281 31 31 0 6.28319}%
    \special{sh 1.000}%
    \special{ar 531 1031 31 31 0 6.28319}%
    \special{sh 1.000}%
    \special{ar 531 781 31 31 0 6.28319}%
    \special{sh 1.000}%
    \special{ar 531 531 31 31 0 6.28319}%
    \special{sh 1.000}%
    \special{ar 781 1281 31 31 0 6.28319}%
    \special{sh 1.000}%
    \special{ar 781 1031 31 31 0 6.28319}%
    \special{sh 1.000}%
    \special{ar 781 781 31 31 0 6.28319}%
    \special{pa 31 1031}%
    \special{pa 31 31}%
    \special{fp}%
    \special{pa 281 1031}%
    \special{pa 281 281}%
    \special{fp}%
    \special{pa 531 1031}%
    \special{pa 531 531}%
    \special{fp}%
    \special{pa 781 1281}%
    \special{pa 781 781}%
    \special{fp}%
    \special{pa 31 281}%
    \special{pa 781 1031}%
    \special{fp}%
    \special{sh 1.000}%
    \special{ar 1531 1031 31 31 0 6.28319}%
    \special{sh 1.000}%
    \special{ar 1531 781 31 31 0 6.28319}%
    \special{sh 1.000}%
    \special{ar 1531 531 31 31 0 6.28319}%
    \special{sh 1.000}%
    \special{ar 1531 281 31 31 0 6.28319}%
    \special{sh 1.000}%
    \special{ar 1531 31 31 31 0 6.28319}%
    \special{sh 1.000}%
    \special{ar 1781 1031 31 31 0 6.28319}%
    \special{sh 1.000}%
    \special{ar 1781 781 31 31 0 6.28319}%
    \special{sh 1.000}%
    \special{ar 1781 531 31 31 0 6.28319}%
    \special{sh 1.000}%
    \special{ar 1781 281 31 31 0 6.28319}%
    \special{sh 1.000}%
    \special{ar 2031 1031 31 31 0 6.28319}%
    \special{sh 1.000}%
    \special{ar 2031 781 31 31 0 6.28319}%
    \special{sh 1.000}%
    \special{ar 2031 531 31 31 0 6.28319}%
    \special{sh 1.000}%
    \special{ar 2281 1281 31 31 0 6.28319}%
    \special{sh 1.000}%
    \special{ar 2281 1031 31 31 0 6.28319}%
    \special{sh 1.000}%
    \special{ar 2281 781 31 31 0 6.28319}%
    \special{pa 1531 1031}%
    \special{pa 1531 31}%
    \special{fp}%
    \special{pa 1781 1031}%
    \special{pa 1781 281}%
    \special{fp}%
    \special{pa 2031 1031}%
    \special{pa 2031 531}%
    \special{fp}%
    \special{pa 2281 1281}%
    \special{pa 2281 781}%
    \special{fp}%
    \special{pa 1531 281}%
    \special{pa 2281 1031}%
    \special{fp}%
    \special{pa 1531 781}%
    \special{pa 1781 1031}%
    \special{fp}%
    \special{pa 1531 531}%
    \special{pa 2031 1031}%
    \special{fp}%
    \special{sh 1.000}%
    \special{ar 3031 1031 31 31 0 6.28319}%
    \special{sh 1.000}%
    \special{ar 3031 781 31 31 0 6.28319}%
    \special{sh 1.000}%
    \special{ar 3031 531 31 31 0 6.28319}%
    \special{sh 1.000}%
    \special{ar 3031 281 31 31 0 6.28319}%
    \special{sh 1.000}%
    \special{ar 3031 31 31 31 0 6.28319}%
    \special{sh 1.000}%
    \special{ar 3281 1031 31 31 0 6.28319}%
    \special{sh 1.000}%
    \special{ar 3281 781 31 31 0 6.28319}%
    \special{sh 1.000}%
    \special{ar 3281 531 31 31 0 6.28319}%
    \special{sh 1.000}%
    \special{ar 3281 281 31 31 0 6.28319}%
    \special{sh 1.000}%
    \special{ar 3531 1031 31 31 0 6.28319}%
    \special{sh 1.000}%
    \special{ar 3531 781 31 31 0 6.28319}%
    \special{sh 1.000}%
    \special{ar 3531 531 31 31 0 6.28319}%
    \special{sh 1.000}%
    \special{ar 3781 1281 31 31 0 6.28319}%
    \special{sh 1.000}%
    \special{ar 3781 1031 31 31 0 6.28319}%
    \special{sh 1.000}%
    \special{ar 3781 781 31 31 0 6.28319}%
    \special{pa 3031 1031}%
    \special{pa 3031 31}%
    \special{fp}%
    \special{pa 3281 1031}%
    \special{pa 3281 281}%
    \special{fp}%
    \special{pa 3531 1031}%
    \special{pa 3531 531}%
    \special{fp}%
    \special{pa 3781 1281}%
    \special{pa 3781 781}%
    \special{fp}%
    \special{pa 3031 281}%
    \special{pa 3781 1031}%
    \special{fp}%
    \special{pa 3031 1031}%
    \special{pa 3281 781}%
    \special{fp}%
    \special{pa 3031 1031}%
    \special{pa 3531 781}%
    \special{fp}%
    \special{pa 3031 781}%
    \special{pa 3281 1031}%
    \special{fp}%
    \special{pa 3031 781}%
    \special{pa 3531 1031}%
    \special{fp}%
    \hbox{\vrule depth1.313in width0pt height 0pt}%
    \kern 3.813in
  }%
}%
\ht\graph=\dp\graph
\dp\graph=0pt

\[\box\graph\]
\caption{$P_4$ and two alternatives.}
\label{F:other}
\end{figure}

We will see (Remark~\ref{R:oursaremin})
that the posets $P_j$ have the minimum cardinality
among all polyunsaturated posets with the same height.
The first two, $P_1$ and $P_2$, are unique;
thus, our $P_2$ is isomorphic to the first of West's posets and to the
Greene-Kleitman example.
However, the rest of the $P_j$'s are not unique.
Figure~\ref{F:other} shows $P_4$ and two other polyunsaturated posets
having the same parameters.

\section{Polyunsaturated Graphs and the Main Result} \label{S:mainres}

We denote the vertex set of a (finite, simple, undirected)
graph $G$ by $V(G)$.
The \emph{independence number} of $G$,
denoted $\alpha(G)$, is the maximum size of an independent set in $G$.
The \emph{clique number} of $G$, denoted $\omega(G)$, is the maximum size
of a clique of $G$.
We refer to~\cite{BoMu76,WesD96} for graph-theoretic terminology not
defined here.

The antichains in a poset $P$ are precisely the independent sets in
its comparability graph $G(P)$;
the chains in $P$ are the cliques in $G(P)$.
This allows us to extend poset properties to graphs.

Given a graph $G$ and a positive integer $k$,
we denote the maximum size of a union of $k$
independent sets in $G$ by $\alpha_k(G)$.
A partition $\cC$ of $V(G)$ into cliques
is \emph{$k$-saturated} if $\alpha_k(G)=m_k(\cC)$.
We say that $G$ is a
\emph{strong Greene-Kleitman graph} (SGK graph) if every induced
subgraph of $G$ has a $k$- and $k+1$-saturated clique
partition, for each positive integer $k$.
An SGK graph $G$ with clique
number $c$ is \emph{polyunsaturated} if $G$ has no $k$-
and $l$-saturated clique partition for any distinct, nonconsecutive
$k,l<c$
(note that every clique partition is $k$-saturated
for each $k\ge c$).

A chain partition $\cC$ of a poset $P$ is
$k$-saturated if and only if $\cC$ is a
$k$-saturated clique partition of $G(P)$.
Thus, $G(P)$ is an SGK graph, by the Greene-Kleitman Theorem.
A poset $P$ is polyunsaturated if and only if $G(P)$ is polyunsaturated.

In this and the following sections, we will be proving
necessary and sufficient conditions for the existence of
polyunsaturated posets.
We will state our results in the more general
graph-theoretic form.
However, for \emph{every}
graph result in the remainder of this paper,
there is an analogous poset result, which can be obtained
by replacing graph-theoretic terminology with the
corresponding poset terminology;
e.g., replace \emph{clique number} with \emph{height},
$\underline\alpha$  with $\underline d$, etc.

The following lemma is due to Greene~\cite[remark after Thm.~3.1]{GreC76}
(see also~\cite[Thm.~4.14]{SakM86}).

\begin{lemma}[Greene 1976] \label{L:ninc}
For every SGK graph $G$, $\Delta\underline\alpha(G)$
is a nonincreasing sequence of positive
integers.\ggcnopf\end{lemma}

This follows from comparing the expressions
$\alpha_k(G)=m_k(\cC)$,
$\alpha_{k+1}(G)=m_{k+1}(\cC)$, and
$\alpha_{k+2}(G)\le m_{k+2}(\cC)$ for a $k$- and $k+1$-saturated
partition $\cC$ and noting
that $\Delta\underline m(\cC)$ is nonincreasing.

Conversely, for each nonincreasing finite sequence
$\underline b$ of positive integers, there is an SGK graph $G$ with
$\Delta\underline\alpha(G)=\underline b$; for example, we can
let $G$ be a disjoint union of cliques of the proper sizes.
Thus we have characterized those sequences that are the
$\Delta\underline\alpha$ sequence of an SGK graph.
Our main result (Thm.~\ref{T:main}) characterizes
those nonincreasing sequences that are
the $\Delta\underline\alpha$ sequence of a polyunsaturated SGK
graph.

The following lemma will be used in this characterization.
A poset version of this lemma was proven
by Greene and Kleitman~\cite[Lemma~3.7]{GrKl76}.

\begin{lemma} \label{L:k+1}
Let $G$ be an SGK graph, and let $k$ be a positive integer.
If $\Delta\alpha_k(G)=\Delta\alpha_{k+1}(G)$,
and $\cC$ is a $k$-saturated clique
partition, then $\cC$ is $k+1$-saturated.
\end{lemma}

\begin{proof}
For all $i$, $\alpha_i(G)\le m_i(\cC)$.
Thus, since
$\cC$ is $k$-saturated, we have
\begin{enumerate}
\item$\alpha_{k-1}(G)\le m_{k-1}(\cC)$, \label{I:ninc-k-1}
\item$\alpha_k(G)=m_k(\cC)$, and \label{I:ninc-k}
\item$\alpha_{k+1}(G)\le m_{k+1}(\cC)$. \label{I:ninc-k+1}
\end{enumerate}
By the definition of $\underline m(\cC)$,
$\Delta m_i(\cC)$ is the number of
cliques in $\cC$ with at least $i$ vertices.
Thus, $\Delta\underline m$ is nonincreasing, and we have
\begin{alignat*}{2}
\Delta\alpha_{k+1}(G)&\le\Delta m_{k+1}(\cC)&
&\quad\text{by subtracting~\ref{I:ninc-k}
  from~\ref{I:ninc-k+1}}\\
&\le\Delta m_k(\cC)&
&\quad\text{since $\Delta\underline m$ is nonincreasing}\\
&\le\Delta\alpha_k(G)&
&\quad\text{by subtracting~\ref{I:ninc-k-1}
  from~\ref{I:ninc-k}}\\
&=\Delta\alpha_{k+1}(G).&
&
\end{alignat*}
Thus,
$\Delta\alpha_{k+1}(G)=\Delta m_{k+1}(\cC)$.
Adding corresponding sides to~\ref{I:ninc-k}, we see that
$\alpha_{k+1}(G)=m_{k+1}(\cC)$, and so
$\cC$ is $k+1$-saturated.\ggcendpf\end{proof}

Our main result is the following.

\begin{theorem} \label{T:main}
Let $c$ be a positive integer and let
$\underline b=\left(b_1,b_2,\dotsc,b_c\right)$ be a nonincreasing
sequence of positive integers.
There exists a polyunsaturated SGK graph $G$
with clique number $c$ and $\Delta\underline\alpha(G)=\underline b$
if and only if $b_2>b_3>\dotsb>b_{c-1}$.
Moreover, we may require both $G$ and $\oG$ to be
comparability graphs.\end{theorem}

\begin{proof}
($\Longrightarrow$)
Let $G$ be a polyunsaturated SGK graph with clique number $c$ and
$\Delta\underline\alpha(G)=\underline b$.
We show that $b_2>b_3>\dotsb>b_{c-1}$.

By Lemma~\ref{L:ninc} the sequence is nonincreasing.
If $\Delta\alpha_k(G)=\Delta\alpha_{k+1}(G)$ with $2\le k\le c-2$,
then Lemma~\ref{L:k+1} implies that every $k-1$- and $k$-saturated
clique partition is also $k+1$-saturated.
The hypothesis that $G$ is
polyunsaturated forbids this, so $b_k>b_{k+1}$ throughout.

($\Longleftarrow$)
Let $\underline b=\left(b_1,\dotsc,b_c\right)$ be a nonincreasing
sequence of positive integers with $b_2>b_3>\dotsb>b_{c-1}$.
It suffices to show that there exists a
polyunsaturated poset $R$ of dimension at most $2$ and height $c$
such that $\Delta\underline d(R)=\underline b$;
$G(R)$ will be the required graph.
If $c<3$, then
every poset with height $c$ is polyunsaturated, and
we may let $R$ be a disjoint union of chains of the appropriate sizes.
Thus, we may assume $c\ge3$.

Since $b_{c-1}\ge b_c\ge1$, we have $b_{c-2}\ge2$;
for $1<i<c$, we have $b_i\ge c-i$.
Also, $b_1\ge b_2\ge c-2$.
Thus,
\[
\sum_{i=1}^c b_i
  \ge c-2
  +\left[\sum_{i=2}^{c-1} c-i\right]
  +1
=\binom{c}{2}.
\]

We proceed by induction on $\sum b_i$, with the base case being
$\sum b_i=\binom{c}{2}$.
If $\sum b_i=\binom{c}{2}$, then the elements of $\underline b$
must all equal the lower bounds found above.
Let $R=P_{c-2}$.
By Statement~\ref{I:mineprops-delta} of Lemma~\ref{L:mineprops},
we have $\Delta\underline d(R)=\underline b$.
By Lemma~\ref{L:minepp}, $R$ is polyunsaturated.
By Statement~\ref{I:mineprops-dim2} of Lemma~\ref{L:mineprops},
$R$ has dimension at most $2$, and so $R$ is the required poset.

Now suppose $\sum b_i>\binom{c}{2}$.
Let $t$ be maximal such
that $b_t$ exceeds the lower bound found above.
We define a new sequence
$\underline b'=\left(b'_1,\dotsc,b'_c\right)$ as follows:
\[
b'_i=
  \cases
    b_i-1,&i\le t;\\
    b_i,&i>t.
  \endcases
\]
This new sequence satisfies
$b'_1\ge b'_2>b'_3>\dotsb>b'_{c-1}\ge b'_c\ge1$.
By the induction hypothesis, there is a polyunsaturated
poset $R'$ of dimension at most $2$ and height
$c$ such that $\Delta\underline d(R')=\underline b'$.
Let $R$ be the disjoint union of $R'$ and a chain of $t$ elements.
Being a disjoint union of posets of dimension at most $2$, $R$ has
dimension at most $2$.
Since $t\le c$, $R$ has height $c$.
If $k$ and $l$ are distinct, nonconsecutive, positive integers
less than $c$, then every $k$- and $l$-saturated chain
partition of $R$ will give such a partition of $R'$.
Since $R'$ is
polyunsaturated, no such partition exists, and so $R$ is
polyunsaturated.
Lastly, $\Delta\underline d(R)=\underline b$, and so $R$ is the required
poset.\ggcendpf\end{proof}

Lemma~\ref{L:ninc} and
Theorem~\ref{T:main} characterize those sequences $\underline b$
with $\underline b=\Delta\underline\alpha(G)$ for \emph{some}
polyunsaturated SGK graph $G$.
However, it is not true that \emph{every} SGK graph $G$ with
$\underline b=\Delta\underline\alpha(G)$ must be polyunsaturated.
Indeed, as noted after Lemma~\ref{L:ninc}, every nonincreasing
sequence $\left(b_1,\dotsc,b_c\right)$ equals
$\Delta\underline\alpha(G)$
for some graph $G$ that is a disjoint union of cliques.
When $c\ge4$, no such graph is polyunsaturated.

Theorem~\ref{T:main} gives necessary and sufficient conditions
for the existence of a polyunsaturated SGK graph $G$ with certain
parameters and then notes that we may require $G$ to be
a comparability graph.
Thus, an analogous result holds for posets.
Translating Theorem~\ref{T:main} into the language of posets,
we obtain the following.

\begin{corollary} \label{C:mainposet}
Let $c$ be a positive integer and let
$\underline b=\left(b_1,b_2,\dotsc,b_c\right)$ be a nonincreasing
sequence of positive integers.
There exists a polyunsaturated poset $P$
with height $c$ and $\Delta\underline d(P)=\underline b$
if and only if $b_2>b_3>\dotsb>b_{c-1}$.
Moreover, we may require $P$ to have dimension
at most $2$.\ggcnopf\end{corollary}

\section{Independence and Clique Numbers} \label{S:icn}

Using Theorem~\ref{T:main},
we can find necessary and sufficient conditions
for the existence of a polyunsaturated SGK graph with prescribed
clique number, independence number, and number of vertices.
As noted earlier, for each result in this section,
an analogous poset result holds.

\begin{corollary} \label{C:3cond}
Let $n$, $c$, $a$ be positive integers, with $c\ge3$.
There exists an $n$-vertex polyunsaturated SGK graph $G$
with clique number $c$ and independence number $a$ if and only if
all of the following conditions hold:
\begin{enumerate}
\item $a\ge c-2$, \label{I:3cond-ac-2}
\item $n\ge a+1+\binom{c-1}{2}$, and \label{I:3cond-nge}
\item $n\le ca+1-\binom{c-1}{2}$. \label{I:3cond-nle}
\end{enumerate}
Moreover, we may require both $G$ and $\oG$ to be comparability graphs.
\end{corollary}

\begin{proof}
By Lemma~\ref{L:ninc} and Theorem~\ref{T:main},
for every sequence $\left(b_1,b_2,\dotsc,b_c\right)$ of positive
integers, there is a polyunsaturated SGK graph $G$ with
$\Delta\underline\alpha(G)=\underline b$ if and only if
\[
b_1\ge b_2>b_3>\dotsb>b_{c-1}\ge b_c,
\]
and this remains true if we require both $G$ and $\oG$
to be comparability graphs.
For every $n$-vertex SGK graph $G$ with clique number $c$
and independence number $a$,
we have
$n=\sum\Delta\underline\alpha(G)$ and $a=\Delta\alpha_1$.
Thus, to prove Corollary~\ref{C:3cond}, it suffices to show that,
for positive integers $n$, $c$, and $a$,
there exists a sequence
$\underline b=\left(b_1,b_2,\dotsc,b_c\right)$ of positive
integers such that
\[
b_1\ge b_2>b_3>\dotsb>b_{c-1}\ge b_c,\quad
  n=\sum_{i=1}^c b_i,\quad\text{and}\quad
  a=b_1
\]
if and only if
conditions~\ref{I:3cond-ac-2}--\ref{I:3cond-nle} hold.

($\Longrightarrow$)
Let $n$, $c$, $a$, and $\underline b$ be as above.
As in the proof of Theorem~\ref{T:main}, since
$b_{c-1}\ge b_c\ge1$, we have $b_i\ge c-i$, for $1<i<c$.
In particular, $a=b_1\ge b_2\ge c-2$.
Thus~\ref{I:3cond-ac-2} holds.
Noting that $b_1\ge a$
(since the two are equal), we can sum all these lower bounds
to obtain~\ref{I:3cond-nge}:
\[
n=\sum_{i=1}^c b_i
  \ge a
  +\left[\sum_{i=2}^{c-1} c-i\right]
  +1
=a+1+\binom{c-1}{2}.
\]

We prove~\ref{I:3cond-nle} similarly, by summing upper bounds.
Since $b_2\le b_1\le a$, we have $b_3\le a-1$, and
we have $b_i\le a-i+2$, for $1<i<c$.
In particular, $b_c\le b_{c-1}\le a-c+3$.
Summing these upper bounds, we obtain
\[
n=\sum_{i=1}^c b_i
  \le a
  +\left[\sum_{i=2}^{c-1}a-i+2\right]
  +a-c+3
=ca+1-\binom{c-1}{2}.
\]

($\Longleftarrow$)
Let $n$, $c$, and $a$ be positive integers satisfying
conditions~\ref{I:3cond-ac-2}--\ref{I:3cond-nle}.
As shown above, if all elements of
$\underline b$ equal their lower bounds---which
requires~\ref{I:3cond-ac-2}---then
$\sum b_i=a+1+\binom{c-1}{2}$.
Similarly,
if all elements of $\underline b$ equal their upper bounds---which
also requires~\ref{I:3cond-ac-2}---then
$\sum b_i=ca+1-\binom{c-1}{2}$.
By~\ref{I:3cond-nge} and~\ref{I:3cond-nle},
$n$ is between these two values, inclusive.
Thus, we may construct a sequence $\underline b$ that
satisfies the required conditions.
We begin by letting each $b_i$ equal
its lower bound.
We increase $b_2$ until either $b_2$ reaches its upper
bound, or $\sum b_i=n$.
Then we increase $b_3$ in a similar
manner, and so on.\ggcendpf\end{proof}

\begin{corollary} \label{C:ac-2}
Let $c$ and $a$ be positive integers.
There exists a polyunsaturated SGK graph $G$ with clique
number $c$ and independence number $a$ if and only if $a\ge c-2$.
Moreover, we may require both $G$ and its
complement $\oG$ to be comparability graphs.
\end{corollary}

\begin{proof}
The cases $c=1$ and $c=2$ are easy.
We assume $c\ge3$.
The necessity follows immediately from Corollary~\ref{C:3cond}.
To show the sufficiency, note that
if $a\ge c-2$, then $a+1+\binom{c-1}{2}\le ca+1-\binom{c-1}{2}$.
Thus, we can choose an $n$ satisfying the conditions of
Corollary~\ref{C:3cond}.\ggcendpf\end{proof}

The following result also follows easily from Corollary~\ref{C:3cond}.

\begin{corollary} \label{C:ncch2}
Let $n$ and $c$ be positive integers, with $c\ge3$.
There exists an $n$-vertex polyunsaturated SGK graph $G$
with clique number $c$ if and only if $n\ge\binom{c}{2}$.
Moreover, we may require both $G$ and $\oG$ to be comparability
graphs.\ggcnopf\end{corollary}

\begin{remark} \label{R:oursaremin}
By the above corollaries, the posets $P_j$ of Definition~\ref{D:mine}
have the minimum width and cardinality over all polyunsaturated posets
with the same height.\ggcendrem\end{remark}

\section{Partitions into Independent Sets} \label{S:indsetpart}

Reversing the roles of cliques and independent sets in our definitions,
we may partition $V(G)$ into independent sets
(such a partition is called a \emph{proper coloring}) and
consider the natural upper bound placed on $\omega_k(G)$, the maximum
number of vertices in a union of $k$ cliques.
We call a proper coloring $\cC$ \emph{$k$-saturated} if this
bound is achieved, that is, if $\omega_k(G)=m_k(\cC)$.
We ask if similar results to those above hold for $k$-saturated
colorings.

A proper coloring $\cC$ of $G$ is $k$-saturated if and only if
$\cC$ is a $k$-saturated clique partition of $\oG$.
Let us call $G$ \emph{co-polyunsaturated} if $\oG$ is
polyunsaturated.
A result of Greene~\cite[Thm.~3.1]{GreC76} states that the complement of
an SGK graph is also an SGK graph.
Thus, for each of our results in Sections~\ref{S:mainres}
and~\ref{S:icn}, there is a dual result.
Some of these are given below.

As in previous sections, for each result below there in an
analogous poset result dealing with partitions of a poset
into antichains.

\begin{corollary} \label{C:main-dual}
Let $a$ be a positive integer and let
$\underline b=\left(b_1,b_2,\dotsc,b_a\right)$ be a nonincreasing
sequence of positive integers.
There exists a co-polyunsaturated SGK graph $G$
with independence number $a$ and
$\Delta\underline\omega(G)=\underline b$
if and only if $b_2>b_3>\dotsb>b_{a-1}$.
Moreover, we may require both $G$ and $\oG$ to be comparability
graphs.\ggcnopf\end{corollary}

\begin{corollary} \label{C:3cond-dual}
Let $n$, $a$, $c$ be positive integers, with $a\ge3$.
There exists an $n$-vertex polyunsaturated SGK graph $G$
with independence number $a$ and clique number $c$ if and only if
all of the following conditions hold:
\begin{enumerate}
\item $c\ge a-2$, \label{I:3cond-dual-ac-2}
\item $n\ge c+1+\binom{a-1}{2}$, and \label{I:3cond-dual-nge}
\item $n\le ac+1-\binom{a-1}{2}$. \label{I:3cond-dual-nle}
\end{enumerate}
Moreover, we may require both $G$ and $\oG$ to be comparability graphs.
\end{corollary}

\begin{corollary} \label{C:ac-2-dual}
Let $a$ and $c$ be positive integers.
There exists a co-polyunsaturated SGK graph $G$ with
independence number $a$ and clique number $c$ if and only if $c\ge a-2$.
Moreover, we may require both $G$ and $\oG$ to be comparability
graphs.\ggcnopf\end{corollary}

\end{document}